\documentclass{article}
\usepackage[a4paper, margin=3cm]{geometry}
\usepackage{amssymb}
\usepackage{graphicx} 
\usepackage{mathtools}
\usepackage{amsthm}
\usepackage{xypic}
\usepackage[colorlinks=true, linkcolor=blue, citecolor=blue, urlcolor=blue]{hyperref}
\usepackage[backend=biber,]{biblatex}
\usepackage{stmaryrd}

\makeatletter
\renewcommand\subsubsection{\@startsection{subsubsection}{3}{\z@}%
  {-18pt plus -4pt minus -4pt}%
  {0pt}%
  {\normalfont\normalsize\bfseries}}
\makeatother

\numberwithin{equation}{section}

\newcommand{\bC}{\mathbb{C}}

\newcommand{\fb}{\mathfrak{b}}

\newcommand{\fg}{\mathfrak{g}}
\newcommand{\fh}{\mathfrak{h}}

\newcommand{\fk}{\mathfrak{k}}

\newcommand{\fn}{\mathfrak{n}}
\newcommand{\fL}{\mathfrak{L}}

\newcommand{\cK}{\mathcal{K}}

\newcommand{\cO}{\mathcal{O}}

\newcommand{\B}{\textnormal{Bun}}

\newcommand{\lpalgp}{\fg(\cO)}
\newcommand{\lpalg}{\fg(\cK)}
\newcommand{\lpgrpp}{G(\cO)}
\newcommand{\lpgrp}{G(\cK)}

\title{Hecke modifications of conformal blocks outside the critical level}
\author{Raschid Abedin, Giovanni Felder, Robert Windesheim}
\date{}

\begin{document}

\maketitle

\begin{abstract}
We define Hecke modifications of conformal blocks over affine Lie algebras at non-critical level by using the Hecke modifications of the underlying $G$-bundles. We show that this procedure is equivalent to the insertion of a twisted vacuum module at an additional marked point and provide an explicit description using coordinate transformations.
\end{abstract}

\section{Introduction} 

Hecke actions play a fundamental role in various incarnations of Langlands correspondences. 
For instance, consider the geometric Langlands correspondence for a simple complex algebraic group \(G\) as proposed by Beilinson-Drinfeld in \cite{BeilinsonDrinfeld} and completed recently in \cite{GLC1,GLC2,GLC3,GLC4,GLC5}. This correspondence states an equivalence between the derived category of critically twisted \(\mathcal{D}\)-modules on the moduli space of \(G\)-bundles \(\B_G^X\) on a complex smooth projective curve \(X\) and the derived category of certain quasi-coherent sheaves on the moduli space of local systems on \(X\) over the Langlands dual group \({}^LG\) of \(G\). 
It is now central that this correspondence is compatible with Hecke actions, which are described by functors called Hecke
modifications, via a specific eigenvalue property. 

The original approach to the geometric Langlands correspondence by \cite{BeilinsonDrinfeld} is closely related to the quantization of the integrable system on \(T^*\B_G^X\) introduced by Hitchin in \cite{hitchin_systems}. 
More recently, an analytic version of the Langlands
program was developed in, e.g., \cite{Teschner2018, EtingofFrenkelKazhdan2022,
  EtingofFrenkelKazhdan2023, EtingofFrenkelKazhdan2024}, which is related to solving the quantized Hitchin system. There, Hecke
modifications define the self-adjoint operators called Hecke operators, which commute with the
Hamiltonians of the quantized Hitchin system.

Both the geometric and analytic Langlands programs can be understood
from the point of view of conformal field theories associated with \(G\) and the representation
theory of the affine Lie algebra \(\widehat{\fg}\) obtained from the Lie algebra \(\fg\) of \(G\). In the geometric case, the sections of
the aforementioned $\mathcal{D}$-modules on \(\B_G^X\) are (holomorphic) conformal blocks of the
Wess--Zumino--Witten model at the critical level, while in the
analytic case, the half-densities are correlation functions combining holomorphic and
antiholomorphic conformal blocks. Here, ``critical level'' means that the value of the central element of the affine Lie algebra \(\widehat{\fg}\) is fixed to be the negative dual Coxeter number of \(G\). For a review of the connections between the geometric Langlands correspondence and conformal field theory, see e.g.\  \cite{Frenkel2007}.

In the aforementioned representation theoretic approach, one usually considers \(\widehat{\fg}\)-modules that are integrable over the positive part \(G(\cO)\) of the loop group.
In settings where one wants to consider more general $\widehat{\fg}$-modules, it can become necessary to consider additional structures on the underlying $G$-bundles at marked points of $X$. For example, if one wants to consider \(\widehat{\fg}\)-modules induced by Verma modules of \(\fg\), one needs to consider $G$-bundles with parabolic structures. This extends \(\B_G^X\) by parameter spaces isomorphic to flag varieties at the marked points.

It would be expected that at least part of the story outlined above extends to non-critical
levels. The main difference is that the center of the universal enveloping algebra of an affine Kac-Moody algebra is non-trivial precisely at the critical level. At the critical level, said center is infinite-dimensional and provides the quantum Hitchin system with sufficient commuting differential operators for integrability; see e.g.\ \cite{frenkel_langlads_correspondance_for_koop_groups}. However, even at non-critical levels, one obtains a non-trivial $\mathcal{D}$-module structure on conformal
blocks by varying marked points on the base curve or the base curve itself, which are both kept fixed in the critical case. Indeed,
the commuting differential operators coming from the quadratic part of
the center of the universal enveloping algebra at the critical level
become projective connections on the sheaf of conformal blocks on the
moduli spaces of curves with marked points at non-critical levels.

In this paper, we examine Hecke modifications for conformal blocks
at non-critical levels with additional structures at marked points of a fixed base curve \(X\). The starting point of our investigation was the
paper by Jeong, Lee, and Nekrasov \cite{jeong_lee_nekrasov}, who
considered the case where \(X = \mathbb{P}^1\), the additional structures are parabolic, and \(G = \textnormal{PGL}_n(\mathbb{C})\). There, they
related Hecke modifications of conformal blocks, understood by Hecke modifying the $G$-bundles underlying the conformal blocks, to insertions of
twisted vacuum modules at an additional marked point of the base curve. Moreover, they described these conformal blocks with $(N+1)$ module insertions using conformal blocks with $N$ insertions and an explicit coordinate transformation. For \(G = \textnormal{PGL}_2(\mathbb{C})\), similar results were obtained by Teschner in \cite{teschner_notes}, in the development of the analytic Langlands correspondence. 

In this paper, we formalize and generalize this approach to Hecke modifications at non-critical levels. In particular, we introduce Hecke modifications of conformal blocks for a very general class of smooth representations of affine Lie algebras at non-critical levels. It is a variation of the usual Hecke modification as defined in \cite{BeilinsonDrinfeld}: instead of a pull-push along the Hecke correspondence, it turns out that the pull-back of twisted \(\mathcal{D}\)-modules in the sense of \cite[Section 1.4]{beilinson_bernstein_jantzen_conjecture} along the Hecke modification map, while keeping the Hecke parameter as an additional structure, is appropriate for us. In Theorem \ref{thm:Hecke_and_insertion}, we prove that this definition is equivalent to the insertion of a twisted vacuum module at an additional marked point, as suggested by \cite{jeong_lee_nekrasov,teschner_notes}. Moreover, we generalize the description of the Hecke modified $(N+1)$-point conformal blocks from $N$-point conformal blocks via coordinate transformations from \cite{jeong_lee_nekrasov,teschner_notes}; see Theorem \ref{thm:N_to_N+1_P^1case}, Example \ref{exp:N_to_N+1_P^1case}, and Theorem \ref{thm:N_to_N+1_general_case}.

This paper is organized as follows. In Section \ref{sec:background}, we review the
background on affine Lie algebras and moduli spaces of
$G$-bundles on smooth projective curves. Twisted coinvariants and conformal
blocks are explained in Section \ref{sec:coinvariants}. In Section \ref{sec:hecke_modification}, we discuss Hecke
modifications of conformal blocks at non-critical levels. Then, our main results are presented in
Section \ref{sec:results}. Finally, in Section \ref{sec:explicit_calculations}, we illustrate our formalism with explicit calculations in the case of the group \(\textnormal{PGL}_2(\mathbb{C})\) with base curve $X = \mathbb P^1$. We recover the formulas from \cite{jeong_lee_nekrasov,teschner_notes} in this case.

\subsection*{Acknowledgments} 
We are grateful to J\"org Teschner for his inspiring explanations and
for sharing his unpublished notes \cite{teschner_notes}.  G.F.\ also thanks the organizers and participants of the Workshop on "Geometry and Representation Theory Associated to $G$-Torsors on Curves" at the Brin Mathematics Research Center, particularly Chiara Damiolini,
Swarnava Mukhopadhay, Johan Martens, Christian Pauly, and Christoph
Schweigert for discussions and comments. Moreover, R.A.\ acknowledges the support by the Deutsche Forschungsgemeinschaft (DFG, German Research Foundation) – SFB 1624 – ``Higher structures, moduli spaces and integrability'' – 506632645.

\section{Background}\label{sec:background}

\subsection{Virasoro and affine Lie algebras}
Let \(\cO = \bC[\![t]\!]\) be the ring of formal power series in the formal variable \(t\) and \(\cK = \bC(\!(t)\!) = \cO[t^{-1}]\) be its field of fractions, i.e.\ the field of complex formal (lower bounded) Laurent power series. Furthermore, let 
\begin{equation}
    D 
\coloneqq \textnormal{Spec}(\cO) \supseteq  D\setminus\{(t)\} = \textnormal{Spec}(\cK) \eqqcolon D^\circ     
\end{equation}
be the formal disk and formal disk without origin respectively. 

For any two \(\mathbb{C}\)-schemes \(V\) and \(S\), let us denote by \(V(S)\) the space of morphisms \(S \to V\). If \(S = \textnormal{Spec}(R)\) is the affine scheme associated to a \(\bC\)-algebra \(R\), we also write \(V(S) = V(R)\).
In particular, if \(V\) is the scheme associated to a finite-dimensional complex vector space, we have
\begin{equation}
    V(S) = V(R) = V \otimes R,   
\end{equation}
where \(\otimes = \otimes_\bC\) is the tensor product over the field of complex numbers throughout this text.

\subsubsection{Virasoro algebra.}\label{subsec:virasoro}\quad 
Consider the Lie algebra
\begin{equation}
    \Gamma(D^\circ,TD^\circ) = \textnormal{Der}(\cK) \coloneqq \cK\partial_t    
\end{equation}
of continuous vector fields on \(D^\circ\) or, equivalently, of continuous derivatives of \(\cK\). Here, ``continuous'' refers to the \((t)\)-adic topology on \(\cK\). It has a natural one-dimensional central extension \(\mathfrak{Vir} = \textnormal{Der}(\cK)  \oplus \bC \mathbf{c}\), which is called {Virasoro algebra}. On the topological basis 
\begin{equation}
    \{L_m \coloneqq -t^{m+1}\partial_t\mid m \in \mathbb{Z}\} \cup \{\mathbf{c}\} \subset \mathfrak{Vir},    
\end{equation}
the Lie bracket is given by \([L_m,\mathbf{c}] = 0\) and
\begin{equation}
        [L_m,L_n] = (m-n)L_{m+n} + \frac{1}{12}(m^3-m)\delta_{m+n,0}\mathbf{c}\,,\qquad m,n\in \mathbb{Z}.
    \end{equation}
Usually, one considers categories of representations of \(\mathfrak{Vir}\) on which \(\mathbf{c}\) acts by multiplication with a fixed number \(c \in \bC\). This value is then referred to as the central charge.

\subsubsection{Affine Kac-Moody algebras.}\label{sec:affine_liealg}\quad Let \(\fg\) be a finite-dimensional simple complex Lie algebra. Consider the Lie algebra \(\lpalg = \fg(D^\circ) = \fg \otimes \cK\) of \(\fg\)-valued regular functions on \(D^\circ\), or equivalently, formal Laurent series with coefficients in \(\fg\). This Lie algebra also has a natural one-dimensional central extension \(\widehat{\fg} = \lpalg \oplus \bC \mathbf{k}\), which is the so-called affine Lie algebra associated with \(\fg\). Its Lie bracket is given on the set of topological generators 
\begin{equation}
    \{a_m \coloneqq a \otimes t^m\mid a \in \fg, m\in \mathbb{Z}\} \cup \{\mathbf{k}\} \subset \widehat{\fg} 
\end{equation}
by \([a_m,\mathbf{k}] = 0\) and 
    \begin{equation}    
        [a_m,b_n] = [a,b]_{m+n} + m\delta_{m+n,0}\kappa(a,b)\mathbf{k}\,,\qquad a,b\in \fg,m,n \in \mathbb{Z}.
\end{equation}
Here, \(\kappa\) is the invariant bilinear form of \(\fg\) normalized by the fact that the square length of the maximal root is 2.

As for \(\mathfrak{Vir}\), one usually considers categories of representations of \(\widehat{\fg}\) on which \(\mathbf{k}\) acts by multiplication with a fixed number \(k \in \bC\). This value is then referred to as level.

\subsubsection{Segal-Sugawara construction.}\label{sec:segalsugawara}\quad The two Lie algebras \(\mathfrak{Vir}\) and \(\widehat{\fg}\) can be related to each other using the so-called {Segal-Sugawara construction}; see e.g.\ \cite[Section 2.5.10]{frenkel2004vertex}. More precisely, consider the completed universal enveloping algebra of level \(k \in \bC\), i.e.\
\begin{equation}\label{eq:completed_universal_enveloping_algerba}
    U_k(\widehat{\fg}) \coloneqq \widehat{U}(\widehat{\fg})/(k - \mathbf{k})\widehat{U}(\widehat{\fg}),
\end{equation}
where \(\widehat{U}(\widehat{\fg})\) is the completion of the universal enveloping algebra \(U(\widehat{\fg})\) of \(\widehat{\fg}\) with respect to the system of ideals generated by \(t^n\lpalgp\), \(n \in \mathbb{N}\).

The Segal-Sugawara construction provides an embedding \(\mathfrak{Vir} \to U_k(\widehat{\fg})\) for \(k \neq -h^\vee\), such that \(\mathbf{c}\) acts via multiplication with 
\begin{equation}
    c_k = \frac{k\, \textnormal{dim}(\fg)}{k + h^\vee},
\end{equation}
where \(h^\vee\) is the dual Coxeter number of \(\fg\).

To construct this embedding, one considers
\begin{equation}
    S(z) \coloneqq \frac{1}{2(k+h^\vee)}\sum_{j = 1}^d :\!I_{j}(z)I^{j}(z)\!:\, =  \sum_{n \in \mathbb{Z}}S_n z^{-n-1} \in U_k(\widehat{\fg})[\![z,z^{-1}]\!].
\end{equation}
Here, we used the following notations:
\begin{itemize}
    \item \(\{I_j\}_{j = 1}^d \subset \fg\) is a basis and \(\{I^j\}_{j = 1}^d \subset \fg^*\) is its dual with respect to
    \(\kappa\);
    
    \item \(a(z) \coloneqq \sum_{n \in \mathbb{Z}} a_nz^{-n-1} = \sum_{n \in \mathbb{Z}} (a \otimes t^n)z^{-n-1}\in \widehat{\fg}[\![z,z^{-1}]\!]\) for any \(a \in \fg\);

    \item \(:- :\) is the normal ordering, defined by
    \begin{equation}
        :\! a(z)b(z)\!: = a(z)_+b(z) + b(z)a(z)_-,
    \end{equation}
    where \(a(z)_+\) and \(a(z)_-\) are the parts of \(a(z)\) with non-negative and negative powers of \(z\) respectively.
\end{itemize}
Now one may verify that
\(L_n \mapsto S_n\) defines a Lie algebra embedding \(\mathfrak{Vir} \to U_k(\widehat{\fg})\). More explicitly, we have
\begin{equation}
    S_n = \frac{1}{2(k+h^\vee)} \sum_{j = 1}^d \left(\sum_{m < 0} I_{j,m}I^{j}_{n-m} + \sum_{m\ge 0}I^{j}_{n-m}I_{j,m}\right).
\end{equation}

\subsubsection{Vertex algebra interpretation.}\label{sec:vertex}\quad
The Segal-Sugawara construction can be understood naturally in the language of vertex algebras. Namely, consider the vacuum representation of \(\widehat{\fg}\) 
\begin{equation}
    V_{\lpgrpp,k} \coloneqq \textnormal{Ind}_{\lpalgp \oplus \bC \mathbf{k}}^{\widehat{\fg}} \bC|0\rangle \cong U(t^{-1}\fg[t^{-1}])|0\rangle,
\end{equation}
of level \(k \in \bC\) and \(\fg(\cO)\) acts trivially on \(|0\rangle\) while \(\mathbf{k}|0\rangle = k|0\rangle\). This representation admits a natural vertex algebra structure \(Y = Y(-,z)\colon V_{\lpgrpp,k} \to \textnormal{End}(V_{\lpgrpp,k})[\![z,z^{-1}]\!]\) determined by \(a_{-1} \mapsto a(z)\). Now, \(S(z) = Y(\omega,z)\) holds for 
\begin{equation}
    \omega \coloneqq \frac{1}{2(k+h^\vee)}\sum_{\alpha = 1}^d I_{\alpha,-1}I^{\alpha}_{-1}|0\rangle \in V_{\lpgrpp,k}.    
\end{equation}

\subsection{Loop groups} 
For the remainder of the paper, we fix a connected simple complex algebra group \(G\) and denote its Lie algebra by \(\fg\).

The Lie algebras \(\lpalgp\) and \(\lpalg\) intuitively correspond to the groups \(\lpgrpp\) and \(\lpgrp\) of algebraic maps \(D \to G\) and \(D^\circ \to G\) respectively. In order to make this correspondence rigorous, one has to give these abstract groups algebraic structure. To do so, one first notes that \(\lpgrpp\) and \(\lpalgp\) are indeed complex algebraic group schemes of infinite type, while \(\lpgrp\) and \(\lpalg\) are complex affine ind-algebraic groups; see e.g.\ \cite{BeilinsonDrinfeld,Beauville_Laszlo_conformal,faltings_loop_groups}.

One can now define an isomorphism between \(\lpalgp\) (resp.\ \(\lpalg\)) and left-invariant (resp.\ continuous left-invariant) vector fields on \(\lpgrpp\) (resp.\ \(\lpgrp\)), manifesting the algebraic group and Lie algebra correspondence in this case. Indeed, \(\lpalg\) is generated as a Lie algebra from nilpotent elements, and for these, we can define left-invariant derivations in the standard way using the exponential map.

\subsubsection{Exponentiation of pro-nilpotent Lie algebra actions.}\label{sec:exponentiation}\quad
The Lie algebras 
\begin{equation}
    \textnormal{Der}_+(\cO) \coloneqq t^2\textnormal{Der}(\cO) \quad\textnormal{ and }\quad\fg_+(\cO) \coloneqq t\lpalgp    
\end{equation}
are pro-nilpotent Lie algebras, i.e.\ they are projective limits of the nilpotent Lie algebras obtained by truncating powers of \(t\).

Every pro-nilpotent Lie algebra \(\mathfrak{u}\) can be integrated to a pro-unipotent algebraic group \(U\). Indeed, one can simply consider the same base space with multiplication defined by the Baker-Campbell-Hausdorff series; see e.g.\ \cite[Chapter IV]{kumar2012kac}. For \(\textnormal{Der}_+(\cO)\) and \(\fg_+(\cO)\) the associated groups are \(\textnormal{Aut}_+(\cO)\), the group of coordinate transformations of \(\cO\) of the form
\begin{equation}
    t \mapsto t + \sum_{n = 2}^\infty a_nt^n \in \cO\,\qquad \{a_n\}_{n = 2}^\infty \subset \bC,
\end{equation}
and \(G_+(\cO)\), the group of elements \(g \in \lpgrpp\) such that \(g(0) = e\), respectively.

As a consequence, if we have a \(\widehat{\fg}\)-module \(M\) such that a subalgebra \(\mathfrak{u} \subseteq \fg_+(D)\) acts locally nilpotent on \(M\), i.e.\ such that for all \(x \in \mathfrak{u}\) and \(m \in M\) there exists an \(N\in \mathbb{N}\) such that \(x^N m = 0\), the action of \(\mathfrak{u}\) integrates to an action of the group \(U\) associated to \(\mathfrak{u}\). If \(\mathfrak{u}\) is a module over a subalgebra of \(\textnormal{Der}_+(\cO)\) via the Segal-Sugawara construction, then this action also integrates to \(\textnormal{Aut}_+(\cO) \ltimes U\). We will use this to induce the actions of appropriate subgroups of \(\textnormal{Aut}(\cO) \ltimes \lpgrpp\) on representations of \(\widehat{\fg}\).

\subsubsection{Adjoint action of \(\lpgrp\) on \(\widehat{\fg}\).}\label{sec:adjoint_action_on_hatg}\quad
The adjoint action \(\textnormal{Ad}_{\lpalg}\) of \(\lpgrp\) on \(\lpalg\) can be extended to \(\widehat{\fg}\). Indeed, for any \(g \in \lpgrp\) and \(a_m = a \otimes t^m \in \fg(\cK)\), where \(a \in \fg\) and \(m \in \mathbb{Z}\), we can define
\begin{equation}\label{eq:adjoint_action_on_hatg}
    \textnormal{Ad}_{\widehat{\fg}}(g)(a_m) = \textnormal{Ad}_{\lpalg}(g)a_m + \textnormal{res}_0\kappa(g^{-1}\partial_tg,a_m)\mathbf{k}.
\end{equation}
If we additionally put \(\textnormal{Ad}_{\widehat{\fg}}(g) \mathbf{k} = \mathbf{k}\), 
this provides a well-defined Lie algebra automorphism of \(\widehat{\fg}\) for any \(g \in \lpgrp\).

\subsubsection{Central extension \(\widehat{G}\) of \(\lpgrp\).}\label{subsec:widehat_G}\quad Similar to \(\lpalg\), \(\lpgrp\) admits a canonical central extension by the multiplicative group \(\bC^\times\). This ind-affine algebraic group will be denoted by \(\widehat{G}\) and is referred to as the affine (Kac-Moody) group associated with \(G\); see \cite{kumar2012kac} for details. An explicit algebraic construction for \(G = \textnormal{SL}_n(\bC)\) can be found in \cite{Beauville_Laszlo_conformal} and a more analytic construction for general \(G\) in the setting of conformal field theory can be found in e.g. \cite{felder_kzb}.

\subsection{Moduli spaces of \(G\)-bundles and configuration spaces}
Let \(X\) be a smooth complex irreducible projective curve and recall that \(G\) is a connected simple complex algebraic group with Lie algebra \(\fg\).

\subsubsection{\(G\)-bundles.}\label{sec:gbundles}\quad
A \(G\)-bundle \(P\) over \(X\) is an algebraic variety equipped with a right \(G\)-action and a morphism \(\pi_P \colon P \to X\) such that:
\begin{itemize}
    \item \(\pi_P(p\cdot g) = \pi_P(p)\) holds for all \(p \in P\) and \(g \in G\);
    \item \(X\) has an \'etale covering \(U \to X\) such that \(P \times_X U \cong G \times U\) as \(G\)-varieties. 
\end{itemize}
We denote the moduli stack of \(G\)-bundles on \(X\) by \(\B^X_G\).

\subsubsection{Bundles of coordinates and trivializations.}\label{subsec:bundles_of_coordinates}\quad Every \(G\)-bundle \(P\) on \(X\) admits a formal trivialization \(P|_{D_z} \cong G \times D_z\) over the formal disk 
\begin{equation}
    D_z \coloneqq \textnormal{Spec}(\cO_z)\,,\quad \textnormal{ where } \quad  \cO_z \coloneqq \widehat{\mathcal{O}}_{X,z} 
\end{equation}
is the local completion of the sheaf of regular functions \(\cO_X\) on \(X\) at \(z\). This gives rise to the \(\lpgrpp\)-bundle \(P^\sim\) of formal trivializations of \(P\) on the space \(X^\sim\) of formal coordinates of \(X\). 

More precisely, \(X^\sim\) is the \(\textnormal{Aut}(\cO)\)-bundle over \(X\) consisting of pairs \((z,s)\) of points \(z \in X\) equipped with a formal coordinate \(s\) on \(D_z\). Furthermore, \(P^\sim\) is the \(\lpgrpp\)-bundle over \(X^\sim\) consisting of triples \(((z,s),\sigma)\), where \((z,s)\in X^\sim\) and \(\sigma \colon P|_{D_z} \to D_z \times G\) is a formal trivialization. We can identify \(\sigma\) with a section \(\Gamma(D_z,P)\) of \(P\) over \(D_z\) and the choice of coordinate \(s\) provides an identification \(D_z \cong D\) and consequently an action of \(\lpgrpp\) on the set of these sections.

\subsubsection{Additional structures at marked points.}\label{sec:additional_structures}\quad
Let us call a subgroup \(K \subseteq G(\cO)\) \emph{nice} if it is a subgroup scheme that is stable under the \(\textnormal{Aut}(\cO)\)-action on \(G(\cO)\) defined by formal coordinate transformations and if the quotient \(G(\cO)/K\) exists as a scheme. 

For a nice subgroup \(K\subseteq G(D)\), a \(K\)-structure of a \(G\)-bundle \(P\) at a point \(z \in X\) is now a choice of an element \(P^\sim|_z/K_z\), where 
\begin{equation}\label{eq:twisted_K}
    K_z \coloneqq D_z \times_{\textnormal{Aut}(\cO)} K. 
\end{equation}
Let us list some important examples of \(K\)-structures for different \(K\):
\begin{itemize}
    \item If \(K = \{e\}\), a \(K\)-structure at \(z \in X\) is the same as a choice of formal trivialization of \(P\) at \(z\);
    
    \item If \(K = \lpgrpp\), there is no additional structure at \(z\);
    
    \item If \(K = G_+(\cO) \eqqcolon G_1(\cO)\), a \(K\)-structure at \(z\) is the same as a trivialization of \(P\) at \(z\). More generally, if \(K = G_n(\cO) \coloneqq \textnormal{Ker}(G(\cO) \to G(\cO/t^n\cO))\), a \(K\)-structure is a trivialization at \(z\) of \(n\)-th order;
    
    \item If \(K = \mathbf{B}\) is the subgroup of elements \(g \in \lpgrpp\) such that \(g(0) \in B\) for a Borel subgroup \(B \subset G\), a \(K\)-structure at \(z\) is the same as a parabolic structure at \(z\).
\end{itemize}
We can add additional structures with respect to different groups at several points of \(X\) at the same time, as long as these points do not coincide.

\subsubsection{Loop group uniformization.}\label{sec:uniformization}\quad
For any non-empty finite set of points 
\begin{equation}
    \underline{z} = (z_i)_{i = 1}^N \in C^X_N \coloneqq \{(z_i)_{i = 1}^N \in X^N\mid z_i \neq z_j \textnormal{ for all }1\le i\neq j\le N\}    
\end{equation}
such that \(z_i \neq z_j\), any \(G\)-bundle \(P\) is trivializable on the formal neighborhood 
\begin{equation}\label{eq:def_D_underlinex}
    D_{\underline{z}} \coloneqq \bigsqcup_{i = 1}^N D_{z_i} = \textnormal{Spec}(\cO_{\underline{z}})\,,\quad \textnormal{ where }\quad \cO_{\underline{z}} \coloneqq \bigoplus_{i = 1}^N\cO_{z_i},
\end{equation}
of \(\{z_1,\dots,z_N\} \subseteq X\) and on the affine curve
\begin{equation}\label{eq:def_Xcirc}
    X^\circ_{\underline{z}} \coloneqq X\setminus\{z_1,\dots,z_N\}.    
\end{equation}
Therefore, it is isomorphic to the \(G\)-bundle \({P}_g\) given by gluing the trivial bundles on \(X^\circ_{\underline{z}}\) and \(D_{\underline{z}}\) using an appropriate \(g \in G(\cK_{\underline{z}})\) as transition map on
\begin{equation}
    D^\circ_{\underline{z}} \coloneqq D_{\underline{z}} \cap X^\circ_{\underline{z}} = D_{\underline{z}} \setminus \{z_1,\dots,z_N\} = \textnormal{Spec}(\cK_{\underline{z}}).
\end{equation}
Here, if \(\cK_z\) denotes the field of fractions of \(\cO_z\) or, equivalently, the functions on \(D_z^\circ\), we wrote:
\begin{equation}
    \cK_{\underline{z}} = \bigoplus_{i = 1}^N \cK_{z_i}.
\end{equation}
In other words, \(P \cong P_g\) is obtained via the following push-out diagram:
\begin{equation}
    \xymatrix{& G \times X^\circ_{\underline{z}} \ar[rd] &\\ G \times D^\circ_{\underline{z}} \ar[ru]^{\subseteq} \ar[rd]_{g} && P_g \\ & G \times D_{\underline{z}} \ar[ru]&}.
\end{equation}
Clearly, changing the trivializations on \(X^\circ_{\underline{z}}\) and \(D_{\underline{z}}\), which amounts to left multiplication by \(G(X^\circ_{\underline{z}})\) and right multiplication by \(G(D_{\underline{z}})\), provides isomorphic bundles. This results in a set-theoretic identification of the set of \(G\)-bundles with the double quotient \(G(X^\circ_{\underline{z}})\!\setminus\!  G(\cK_{\underline{z}})/G(\cO_{\underline{z}})\). 

This identification can be made algebro-geometrically in the language of stacks in order to obtain the loop group uniformization. In other words, we have an isomorphism
\begin{equation}\label{eq:uniformization}
    \textnormal{Bun}_G^X \cong G(X^\circ_{\underline{z}})\!\setminus\!  G(\cK_{\underline{z}})/G(\cO_{\underline{z}})
\end{equation}
of stacks for every \(\underline{z} \in C^X_N\), where \(\B_G^X\) is the moduli space of \(G\)-bundles on \(X\); see \cite{drinfeld_simpson}. Here, we note that \(G(X_{\underline{z}}^\circ)\!\setminus\!G(\cK_{\underline{z}})\) is actually a scheme (of infinite type), so the right hand side of \eqref{eq:uniformization} refers to the quotient stack of this scheme by the affine infinite-type group scheme \(G(\cO_{\underline{z}})\).

\subsubsection{Uniformization with additional structures.}\label{sec:uniformization_additional_str}\quad One can also understand the moduli stack of \(G\)-bundles with additional structures given by a family \(\underline{K} = (K_1,\dots,K_N)\) of nice subgroups of \(\lpgrpp\) using uniformization. Namely, writing \(K_{\underline{z}} = K_{1,z_1} \times \dots \times K_{N,z_N}\), using the notation of \eqref{eq:twisted_K},
we have an isomorphism of stacks
\begin{equation}
    \textnormal{Bun}_{G,\underline{K}}^{X,\underline{z}} \cong G(X^\circ_{\underline{z}})\!\setminus\!  G(\cK_{\underline{z}})/K_{\underline{z}},
\end{equation}
where \(\textnormal{Bun}_{G,\underline{K}}^{X,\underline{z}}\) is the moduli stack of \(G\)-bundles with \(K_i\)-structures at \(z_i\).

\subsubsection{Important moduli spaces.}\label{subsec:important_moduli}\quad
Let us introduce notations for the remaining moduli spaces we consider in the work.
The largest moduli space we consider is the space \(\textnormal{Bun}^{X,\sim}_{G,N}\) of \(G\)-bundles on \(X\) with formal trivializations at \(N\) marked points. Changing the formal trivializations provides this space with a \(G(\cO)^N\)-torsor structure over
\begin{equation}
    C^{X,\sim}_N \coloneqq X^{\sim,N}|_{C^X_N} =  \{(z_i,s_i)_{i = 1}^N \in X^{\sim,N}\mid z_i \neq z_j, 1\le i\neq j\le N\},
\end{equation}
where \(X^\sim\) was defined in Section \ref{subsec:bundles_of_coordinates}, which is itself a \(\textnormal{Aut}(\cO)^N\)-torsor over the open subset of \(C^X_N\).
Let us note that \(\textnormal{Bun}^{X,\sim}_{G,N}\) is a scheme (of infinite type) whose fiber over \((z_i,s_i)_{i = 1}^N \in C^{X,\sim}_N\) is naturally identified with \(G(X^\circ_{\underline{z}})\!\setminus\!G(\cK_{\underline{z}})\).
The other moduli spaces we consider are essentially quotients of \(\textnormal{Bun}^{X,\sim}_{G,N}\) and \(C^{X,\sim}_N \). 

We have an \(\textnormal{Aut}_+(\cO)\)-action on \(X^\sim\), so that we can pass to 
\begin{equation}
    \begin{split}
         T^\times C^X_N \coloneqq T^\times X^{N}|_{C^X_N} &= \{\{(z_i,v_i)\}_{i = 1}^N \in T^\times X^{N}\mid z_i \neq z_j, 1\le i\neq j\le N\} 
        \\&\cong C^{X,\sim}_N/\textnormal{Aut}_+(\cO)^N.
    \end{split}
\end{equation}
Here, \(T^\times X \subseteq TX\) is the \(\bC^\times\)-bundle of non-vanishing tangent vectors of \(X\), which is isomorphic to \(X^\sim/\textnormal{Aut}_+(\cO)\).

Let \(\underline{K} = (K_1,\dots,K_N)\) be a family of nice subgroups of \(\lpgrpp\) in the sense of Section \ref{sec:additional_structures}. We write \(\textnormal{Bun}_{G,\underline{K}}^{X}\) for the moduli space over \(T^\times C^X_N\) of \(G\)-bundles with \(K_i\)-structures at the \(i\)-th marked point of \(X\). Then 
\begin{equation}\label{eq:add_str_as_quotient}
    \textnormal{Bun}_{G,\underline{K}}^{X} = \textnormal{Bun}_{G,N}^{X,\sim}/\Pi_{i = 1}^N(\textnormal{Aut}_+(\cO) \ltimes K_i)  
\end{equation}
by virtue of Section \ref{sec:uniformization_additional_str}.
If we fix a \(G\)-bundle \(P\), we write $C^P_{\underline{K}}$ for the moduli space of \(\underline{K}\)-structures on \(P\) at \(N\) marked points of \(X\), which consists of the points 
\begin{equation}\label{eq:def_config_P}
    (p_1,\dots,p_N) \in \prod_{i = 1}^N P^\sim/(\textnormal{Aut}_+(\cO) \ltimes K_i) 
\end{equation}
satisfying \(
    (\pi_{P^\sim}(p_1),\dots,\pi_{P^\sim}(p_N)) \in C_N^{X,\sim}\). We can consider \(C^P_{\underline{K}}\) as a subset of \(\textnormal{Bun}_{G,\underline{K}}^X\). It is the fiber over \(P\) under the canonical projection \(\textnormal{Bun}_{G,\underline{K}}^X \to \textnormal{Bun}_G^X\) that forgets about additional structures.

\subsubsection{Line bundles on the moduli spaces.}\quad \label{sec:line_bundles}
The central extension \(\widehat{G}\) of \(G(\cK)\) from Section \ref{subsec:widehat_G} gives rise to a \(\bC^\times\)-bundle \(G(X_{\underline{z}})\!\setminus\!\widehat{G} \to G(X_{\underline{z}})\!\setminus\!G(\cK_{\underline{z}})\). By varying \(\underline{z} \in X^N\), we obtain a \(\bC^\times\)-bundle \(\mathcal{L}^\sim_N\) on \(\textnormal{Bun}_{G,N}^{X,\sim}\). This \(\bC^\times\)-bundle descends to a \(\bC^\times\)-bundle \(\mathcal{L}_{\underline{K}}\) on \(\textnormal{Bun}_{G,\underline{K}}^{X}\).

Let us note that the natural right \(G(\cK)^N\)-action on \(\textnormal{Bun}_{G,N}^{X,\sim}\) extends to a \(\widehat{G}^N\)-action on \(\mathcal{L}^\sim_N\). Therefore, we obtain a \(\widehat{\fg}^N\)-action on \(\mathcal{L}^\sim_N\).
In this work, the \(\mathcal{D}\)-module structures will be twisted using this \(\bC^\times\)-bundle \(\mathcal{L}_N^\sim\) and the \(\widehat{\fg}\)-action on it. 

\section{Coinvariants and conformal blocks}\label{sec:coinvariants}
Recall that we fixed a smooth complex projective curve \(X\) and a simple complex connected algebraic group \(G\) with Lie algebra \(\fg\).

\subsection{Vector bundles associated to \(G\)-bundles}
In this section, we give a brief survey on how we can use \(G\)-bundles and representations of \(G\) in order to construct vector bundles, essentially following \cite{frenkel2004vertex}, but generalizing the approach given there to allow for more general representations of \(\widehat{\fg}\). 

\subsubsection{Vector bundles associated to \(G\)-bundles and \(G\)-modules.}\label{sec:associated_vbundles_Gmod}\quad
Essentially, we will use an infinite-dimensional version of the following well-known construction of vector bundles from \(G\)-modules and \(G\)-bundles.
Let \(V\) be a representation of \(G\) and \(P\) be a \(G\)-bundle. Then 
\begin{equation}
    V^P \coloneqq P \times_G V = (P \times V)/G,    
\end{equation}
where \(G\) acts via \(g \cdot (p,v) = (p\cdot g^{-1},g\cdot v)\), is a vector bundle over \(X\). 

This construction can be generalized if \(V\) is not necessarily a \(G\)-module but only a \(K\)-module for a closed subgroup \(K \subseteq G\). In this case, one can consider
\begin{equation}
    V^P_K \coloneqq P \times_K V \to P/K\,,\qquad [p,v] \mapsto [p].  
\end{equation}
Observe that \(P/K \cong X \times G/K\) if and only if \(P\) admits a reduction to a \(K\)-bundle. In particular, for \(K = G\) we have \(P/G \cong X\) and we recover \(V^P_G = V^P\).

For example, if \(K = B\) is a Borel subgroup, every \(G\)-bundle admits a reduction to a \(B\)-bundle; see \cite{drinfeld_simpson}. Therefore, we have a vector bundle \(V^P_B \to X \times G/B\) for any \(B\)-module \(V\).
In particular, this holds for Verma modules.

\subsubsection{Vector bundles associated to \(G\)-bundles and \(\widehat{\fg}\)-modules I.}\label{sec:vbundles_hatg_modules_I}\quad
Consider a Lie algebra \(\mathfrak{L}\) with a Lie subalgebra \(\mathfrak{K}\) that integrates to an algebraic group \(K\). Recall that \(M\) is called a Harish-Chandra \((\fL,K)\)-module, or simply \((\fL,K)\)-module, if 
\begin{itemize}
    \item[(H1)] \(M\) is an \(\fL\)-module and \(K\)-module;
    \item[(H2)] The actions of \(\mathfrak{K}\) on \(M\) induced by \(\fL\) and \(K\) coincide.
\end{itemize}

Let us call a \(\widehat{\fg}\)-module \(M\) of level \(k \neq -h^\vee\), i.e.\ on which \(\mathbf{k}\in\widehat{\fg}\) acts by multiplication with \(k \in \bC\), \emph{admissible} if:
\begin{itemize}
    \item[(A1)] \(M\) is smooth, i.e.\ \(t^m\lpalgp \subseteq \widehat{\fg}\) acts by zero for sufficiently large \(m \in \mathbb{N}\);

    \item[(A2)] The action of \(\widehat{\fg}\) on \(M\) is ind-algebraic, i.e.\ \(M = \bigcup_{n \in \mathbb{N}} M_n\) for a growing series of finite-dimensional vector subspaces \(M_n \subseteq M\) and for all \(n,m \in \mathbb{N}\) there exists an \(\ell \in \mathbb{N}\) such that \(t^{-n}\lpalgp \cdot  M_m \subseteq M_{\ell}\).

    \item[(A3)] The action of \(\textnormal{Der}_+(\cO)\subseteq \mathfrak{Vir}\) on \(M\) induced from the Segal-Sugawara construction is locally nilpotent.
\end{itemize}
Let us note that smooth \(\widehat{\fg}\)-modules of level \(k\) are precisely \(U_k(\widehat{\fg})\)-modules. Furthermore, (A3) induces an \(\textnormal{Aut}_+(\cO)\)-action on \(M\) via exponentiation; see Section \ref{sec:exponentiation}.

To any admissible
\((\widehat{\fg},\lpgrpp)\)-module \(M\) of level \(k\) and any \(G\)-bundle \(P\), we can associate
\begin{equation}\label{eq:def_MP}
    \begin{split}
    \Delta^{P,\sim}(M) \coloneqq P^\sim \times_{\lpgrpp} ({X}^\sim \times_{\textnormal{Aut}_+(\cO)} M) = P^\sim \times_{\textnormal{Aut}_+(\cO) \ltimes \lpgrpp} M,
    \end{split}
\end{equation}
where \(X^\sim\) and \(P^\sim\) mean the same as in Section \ref{subsec:bundles_of_coordinates}. This defines a (potentially infinite-dimensional) vector bundle
\begin{equation}
    \Delta^{P,\sim}(M) \to T^\times X.
\end{equation}
Here, recall that
\(T^\times X \subseteq TX\) denotes the space of pairs \((x,v)\) of points \(x\in X\) equipped with a non-vanishing tangent vector \(v \in T_xX \setminus \{0\}\) and the canonical projection is defined using
\begin{equation}
    P^\sim/(\textnormal{Aut}_+(\cO) \ltimes \lpgrpp) = X^\sim/\textnormal{Aut}_+(\cO) = T^\times X.
\end{equation}
The space \(\Delta^{P,\sim}(M)\) is a (potentially infinite-dimensional) vector bundle whose fibers are all isomorphic to the ind-algebraic vector space \(M\). 

As a side note: if the action of \(L_0 = -t\partial_t\) additionally acts on \(M\) by scalar multiplication with an integer, we would even obtain an action of the whole automorphism group \(\textnormal{Aut}(\cO)\) on \(M\) and we could define a vector bundle \(\Delta^{P,\sim}(M) \to X\) in a similar fashion as above. However, this assumption is not reasonable for the endeavors of this paper, since the eigenvalue of \(L_0\) varies with the level \(k\) which we want to keep arbitrary; see e.g.\ Lemma \ref{lem:segal_sugawara_conjugation} below.

\subsubsection{Vector bundles associated to \(G\)-bundles and \(\widehat{\fg}\)-modules II.}\label{sec:vbundles_hatg_modules_II}\quad
We can consider more general Harish-Chandra modules. Let \(K\subseteq G(\cO)\) be a nice subgroup in the sense of \ref{sec:additional_structures}. We say that a \(\widehat{\fg}\)-module \(M\) is \emph{\(K\)-admissible}, if \(M\) is a \((\widehat{\fg},K)\)-module and \(M\) is admissible as \(\widehat{\fg}\)-module in the sense of Section \ref{sec:vbundles_hatg_modules_I}.
Using the exponentiation from Section \ref{sec:exponentiation}, we see that \(M\) has a natural action of \( \textnormal{Aut}_+(\cO) \ltimes K\).
In the following, we generalize the procedure from Section \ref{sec:vbundles_hatg_modules_I} to \(K\)-admissible \(\widehat{\fg}\)-modules of level \(k\neq -h^\vee\) for any nice subgroup \(K \subseteq \lpgrpp\).

Let us fix a nice subgroup \(K \subseteq G(\cO)\) and a \(K\)-admissible representation \(M\) of level \(k \neq -h^\vee\) for the rest of Section \ref{sec:coinvariants}.
 We can associate the vector bundle 
\begin{equation}
    \Delta^{P,\sim}_K(M) = P^\sim \times_{\textnormal{Aut}_+(\cO) \ltimes K} M \longrightarrow P^\sim/({\textnormal{Aut}_+(\cO) \ltimes K}) = C^P_K
\end{equation}
to any \(G\)-bundle \(P\) on \(X\).
If \(P^\sim\) admits a reduction to \({\textnormal{Aut}(\cO) \ltimes K}\), then 
\begin{equation}
    {P}^\sim/({\textnormal{Aut}_+(\cO) \ltimes K}) \cong T^\times X \times \lpgrpp/K    
\end{equation}
and for \(K = \lpgrpp\) we recover \(\Delta^{P,\sim}_{\lpgrpp}(M) = \Delta^{P,\sim}(M)\).

Let us note that if \(M\) is \(K\)-admissible, it is clearly also \(K'\)-admissible for every nice subgroup \(K'\subseteq K\). One can recognize the additional structure of \(\Delta^{P,\sim}_{K'}(M)\) coming from the integrability of \(M\) over \(K\) through the pull-back diagram
\begin{equation}
    \xymatrix{\Delta^{P,\sim}_{K'}(M)\ar[r]\ar[d]& C^P_{K'}\ar[d]\\ \Delta^{P,\sim}_{K}(M) \ar[r]& C^P_K}.
\end{equation}

\subsubsection{Vector bundles associated to \(G\)-bundles and \(\widehat{\fg}\)-modules III.}\label{sec:vbundles_hatg_modules_III}\quad
We can further generalize the construction of vector bundles from \(K\)-admissible \(\widehat{\fg}\)-modules by allowing the \(G\)-bundles to vary. In particular, recall that \(\textnormal{Bun}_{G,1}^{X,\sim}\) is the space of \(G\)-bundles with formal trivialization at a marked point of \(X\) which comes equipped with a formal coordinate; see Section \ref{subsec:important_moduli}. Then
\begin{equation}
    \Delta^\sim_K(M) = 
    \textnormal{Bun}_{G,1}^{X,\sim} \times_{\textnormal{Aut}_+(\cO) \ltimes K} M \to \textnormal{Bun}_{G,K}^X
\end{equation}
is a vector bundle over the moduli space \(\textnormal{Bun}_{G,K}^X\) of \(G\)-bundles with \(K\)-structures at the marked point. Indeed, we clearly have \(
    \Delta^{P,\sim}_K(M) = \Delta^\sim_K(M)|_{C^P_K}\).

The sheaf of sections of \(\Delta^\sim_K(M)\) can be written as
\begin{equation}
    \pi_{K,*}(\cO_{\textnormal{Bun}_{G,1}^{X,\sim}} \otimes \widehat{\fg})^{\textnormal{Aut}_+(\cO) \ltimes K}.
\end{equation}
Here, \(\pi_K\colon \textnormal{Bun}_{G,1}^{X,\sim} \to \textnormal{Bun}_{G,K}^{X}\) denotes the canonical projection.

\subsection{Coinvariants and conformal blocks}
 
\subsubsection{Definition.}\label{sec:coinvariants_Ia}\quad
Until now, the \(\widehat{\fg}\)-module structure of a \(K\)-admissible representation \(M\) of level \(k\) was only used to integrate the action over the Lie algebra of \(K\) and the action of the \(\textnormal{Der}_+(\cO) \subset \mathfrak{Vir}\). In order to define the proper space of coinvariants and its dual space, the space of conformal blocks, this will change.

To do so, let \(\widehat{\fg}_{\textnormal{out},1} \subset \mathcal{L}_1^\sim \times \widehat{\fg}\) be the kernel of the canonical Lie algebroid map \(\mathcal{L}_{1}^\sim \times \widehat{\fg} \to T\mathcal{L}_{1}^\sim\),
where \(\mathcal{L}_1^\sim\) is the \(\bC^\times\)-bundle on \(\textnormal{Bun}_{G,1}^{X,\sim}\) from Section \ref{sec:line_bundles}. This kernel descends to a Lie algebra bundle on \(\widehat{\fg}_{\textnormal{out},1} \subset \textnormal{Bun}_{G,1}^{X,\sim} \times \widehat{\fg}\) over \(\textnormal{Bun}_{G,1}^{X,\sim}\). More explicitly, the fiber of \(\widehat{\fg}_{\textnormal{out},1}\) over \(((z,s),(P,\sigma)) \in \textnormal{Bun}_{G,1}^{X,\sim}\) is
\begin{equation}
    \widehat{\fg}^P_{\textnormal{out},z}\coloneqq\textnormal{Ad}_{\widehat{\fg}}(g)\fg(X^\circ_z) \subseteq \widehat{\fg},
\end{equation}
where \(g\) is the representative of \((P,\sigma)\) in \(\textnormal{Bun}_{G,1}^{X,\sim}|_{(z,s)} \cong G(X_z^\circ)\!\setminus\!G(D^\circ_z)\). 

The image of the fiberwise action of \(\widehat{\fg}_{\textnormal{out},1}\) on \(\textnormal{Bun}_{G,1}^{X,\sim} \times M\) is easily seen to be stable under the \((\textnormal{Aut}_+(\cO) \ltimes K)\)-action, and we denote the \((\textnormal{Aut}_+(\cO) \ltimes K)\)-invariants under this image by \(\widehat{\fg}_{\textnormal{out},1}\cdot \Delta^\sim_{{K}}(M)\).
The sheaf of coinvariants on the space \(\textnormal{Bun}_{G,K}^X\) is defined by 
\begin{equation}
    \Delta_{{K}}(M) \coloneqq \Delta^\sim_{{K}}(M)/ \,\widehat{\fg}_{\textnormal{out},1}\cdot \Delta^\sim_{{K}}(M).
\end{equation}
For a fixed \(G\)-bundle \(P\), we get a sheaf of coinvariants over \(C^P_K\) defined by 
\begin{equation}
    \Delta_{{K}}^P(M) \coloneqq \Delta^{P,\sim}_{{K}}(M)/ \,\widehat{\fg}^P_{\textnormal{out},1}\cdot \Delta^{P,\sim}_{{K}}(M) = \Delta_K(M)|_{C^P_K},
\end{equation}
where \(\widehat{\fg}_{\textnormal{out},1}^P \coloneqq \widehat{\fg}_{\textnormal{out},1}|_{C^P_K}\).
The sheaves of conformal blocks \(\mathcal{C}_{{K}}(M)\) and \(\mathcal{C}^P_K(M)\) on \(\textnormal{Bun}_{G,K}^X\) and \(C^P_K\) are the respective dual sheaves.

\subsubsection{Multi-point generalizations.}\label{sec:multiple_points}\quad
All the constructions in this section admit multi-point analogs in a straightforward way. For the remainder of Section \ref{sec:coinvariants}, let us fix a family 
\begin{equation}
    \underline{K} = (K_1,\dots,K_N)    
\end{equation}
of nice subgroups of \(G(\cO)\). Let us call a family \(\underline{M} = (M_1,\dots,M_N)\) of \(\widehat{\fg}\)-modules of level \(k \neq -h^\vee\) such that \(M_i\) is \(K_i\)-admissible for \(i \in \{1,\dots,N\}\) simply \(\underline{K}\)-admissible of level \(k\). We fix such a family for the rest of Section \ref{sec:coinvariants}. 

We can consider
\begin{equation}
    \begin{split}
        &\Delta^{\sim}_{\underline{K}}(\underline{M}) = \textnormal{Bun}_{G,N}^{X,\sim} \times_{\prod_{i = 1}^N(\textnormal{Aut}_+(\cO) \ltimes K_i)} (M_1 \otimes \dots \otimes M_N), \\
        &\Delta^{P,\sim}_{\underline{K}}(\underline{M}) = \left(\Delta^{P,\sim}_{K_1}(M_1) \boxtimes \dots \boxtimes \Delta_{K_N}^{P,\sim}(M_N)\right)\Big|_{C_{\underline{K}}^P},
    \end{split}
\end{equation}
and the associated spaces of coinvariants 
\begin{equation}
    \begin{split}
        &\Delta_{\underline{K}}(\underline{M}) \coloneqq \Delta^\sim_{\underline{K}}(\underline{M})/ \,\widehat{\fg}_{\textnormal{out},N}\cdot \Delta^\sim_{\underline{K}}(\underline{M}),\\
        &\Delta_{\underline{K}}^P(\underline{M}) \coloneqq \Delta^{P,\sim}_{\underline{K}}(\underline{M})/ \,\widehat{\fg}^P_{\textnormal{out},N}\cdot \Delta^{P,\sim}_{\underline{K}}(\underline{M}).    
        \end{split}
\end{equation}
Here, \(\widehat{\fg}_{\textnormal{out},N} \coloneqq \textnormal{Ker}\left( \textnormal{Bun}_{G,N}^{X,\sim} \times \widehat{\fg}_N \to T\mathcal{L}_N^\sim\right)\) and \(\widehat{\fg}_{\textnormal{out},N}^P\) is the fiber over a fixed \(G\)-bundle \(P\), while \(\widehat{\fg}_N\) is the one-dimensional central extension of \(\fg(\cK)^{\oplus N}\) in the same way as for \(\fg(\cK)\) in Section \ref{sec:affine_liealg}. We can embed 
\begin{equation}
    \widehat{\fg}_N = \fg(\cK)^{\oplus N} \oplus \bC \mathbf{k} \to \widehat{\fg}^{\oplus N} = \fg(\cK)^{\oplus N} \oplus \bigoplus_{i = 1}^N \bC \mathbf{k}_i  
\end{equation}
by \(\mathbf{k} \mapsto \frac{1}{N}\sum_{i = 1}^N \mathbf{k}_i\).
The multi-point conformal blocks \(\mathcal{C}_{\underline{K}}(\underline{M})\) and \(\mathcal{C}^P_{\underline{K}}(\underline{M})\) are now the dual sheaves of the respective coinvariants.

\subsubsection{Remark: coinvariants vs.\ conformal blocks.}\quad In the following, we will work mainly with the sheaf of coinvariants instead of the sheaf of conformal blocks, since the constructions in the subsequent sections become more transparent when working with coinvariants. Moreover, it has formal advantages, as the sheaf of coinvariants is actually quasi-coherent, while the sheaf of conformal blocks might fail to be due to the infinite-dimensional setting. However, we want to remark that, by using duality, the constructions in the following can be transported to conformal blocks.

\subsubsection{Alternative definition of coinvariants.}\label{sec:alt_def_of_localzation}\quad
Another convenient way of describing coinvariants, or more precisely, its sheaf of sections, is
\begin{equation}\label{eq:localization}
    \Delta_{K}(M) \coloneqq \left(\pi_{K,*}\left(\mathcal{D}_{1,k}^\sim \otimes_{U_k(\widehat{\fg})} M\right)\right)^{\textnormal{Aut}_+(\cO) \ltimes K},
\end{equation}
and its multi-point analog
\begin{equation}\label{eq:localization_multipoints}
    \Delta_{\underline{K}}(\underline{M}) \coloneqq \left(\pi_{\underline{K},*}\left(\mathcal{D}_{N,k}^\sim \otimes_{U_k(\widehat{\fg}_N)} (M_1 \otimes \dots \otimes M_N)\right)\right)^{\prod_{i = 1}^N(\textnormal{Aut}_+(\cO) \ltimes K_i)}.
\end{equation}
Here, \(\mathcal{D}_{N,k}^\sim\) is the sheaf of \(k\)-twisted differential operators on \(\mathcal{L}_N^\sim\) relative to \(C_{N}^{X,\sim}\) and \(\pi_{\underline{K}} \colon \B_{G,N}^{X,\sim} \to \B^X_{G,\underline{K}}\) is the canonical projection.

\subsubsection{Differential operators as coinvariants.}\label{sec:diff_ops_as_coinvariants}\quad
An important class of examples of coinvariants is given by sheaves of twisted differential operators themselves. Indeed, for any nice subgroup \(K \subseteq G(\cO)\) with Lie algebra \(\mathfrak{k} \subseteq \fg(\cO)\), let 
\begin{equation}\label{eq:general_vacuum}
    V_{K,k} \coloneqq \textnormal{Ind}_{\mathfrak{k} \oplus \bC\mathbf{k}}^{\widehat{\fg}}\bC |0\rangle_\mathfrak{k} \cong \left(U_k(\widehat{\fg})/U_k(\widehat{\fg})\mathfrak{k}\right)|0\rangle_\mathfrak{k}
\end{equation}
be the associated vacuum module of level \(k\). Here, \(\mathfrak{k}\) acts trivially on \(|0\rangle_\mathfrak{k}\) and \(\mathbf{k}|0\rangle_\mathfrak{k} = k|0\rangle_\mathfrak{k}\).
Then 
\begin{equation}
    \Delta_K(V_{K,k}) = \mathcal{D}_{K,k}
\end{equation}
is the sheaf of \(k\)-twisted differential operators on \(\textnormal{Bun}_{G,K}^X\) relative to \(T^\times X\);
see \cite{frenkel_benzvi_geometric_realization}. A way to see this is as follows. By definition, the kernel of the Lie algebroid map \(\mathcal{L}_{1}^\sim \times \widehat{\fg} \to T\mathcal{L}_1^\sim\), where \(\mathcal{L}^\sim_1\) is the \(\bC^\times\)-bundle from Section \ref{sec:line_bundles}, is the pull-back of \(\widehat{\fg}_{\textnormal{out},1}\) to \(\mathcal{L}^\sim_1\). On the other hand, taking \((\textnormal{Aut}_+(\cO) \ltimes K)\)-invariants of \(\mathcal{L}_{1}^\sim \times \widehat{\fg}\) and modding out the trivial action of \(\mathcal{L}_1^\sim \times_{\textnormal{Aut}_+(\cO) \ltimes K} \fk\), we obtain the Lie algebroid \(\mathcal{L}_1^\sim \times_{\textnormal{Aut}_+(\cO) \ltimes K} \widehat{\fg}/\fk\). By modding out the \(\widehat{\fg}_{\textnormal{out},1}\)-action as well, we obtain the tangent algebroid of \(\mathcal{L}_{K}\). The universal enveloping algebroid of this Lie algebroid is thus the sheaf of differential operators on \(\mathcal{L}_{K}\). Now \(\Delta_K(V_{K,k})\) is obtained from this algebroid by quantum Hamiltonian reduction with respect to the \(\bC^\times\)-action and moment value \(k\), which is precisely \(\mathcal{D}_{K,k}\).

There is the following multi-point generalization, which can be deduced similarly, for the family \(\underline{K} = (K_1,\dots,K_N)\) of nice subgroups of \(G(\cO)\): 
\begin{equation}
    \Delta_{\underline{K}}(V_{K_1,k},\dots, V_{K_N,k}) = \mathcal{D}_{\underline{K},k}.
\end{equation}
Here, \(\mathcal{D}_{\underline{K},k}\) is the sheaf of \(k\)-twisted differential operators on \(\textnormal{Bun}_{G,\underline{K}}^X\) relative to \(T^\times C^X_N\).
After fixing a \(G\)-bundle \(P\) on \(X\), this can be adjusted as follows: we can describe the sheaf of \(k\)-twisted differential operators on \(C^P_{\underline{K}}\) relative to \(T^\times C^X_N\) as
\begin{equation}
    \Delta^P_{\underline{K}}(V_{K_1,k},\dots, V_{K_N,k}) = \mathcal{D}^P_{\underline{K},k}.
\end{equation}

\subsubsection{Connection on coinvariants.}\label{sec:connection}\quad
The sheaves associated with \(\widehat{\fg}\)-modules above admit a natural twisted connection coming from the unused \(\widehat{\fg}\) and \(\mathfrak{Vir}\)-actions. More precisely, \(\Delta_{\underline{K}}(\underline{M})\) is naturally a twisted \(\mathcal{D}\)-module over \(\B_{G,\underline{K}}^X\).
According to \cite{frenkel_benzvi_geometric_realization}, we can split the twisted differential operators over \(\B_{G,\underline{K}}^X\) into its \(T^\times C_N^X\)-component, along which we obtain the KZB connection for a fixed curve, and the \(G\)-bundle component, which is described precisely by a \(\mathcal{D}_{\underline{K},k}\)-action, where \(\mathcal{D}_{\underline{K},k}\) was defined in Section \ref{sec:diff_ops_as_coinvariants}. 

Following \cite{frenkel_benzvi_geometric_realization}, the \(\mathcal{D}_{\underline{K},k}\)-action on \(\Delta_{\underline{K}}(\underline{M})\) can be described using Section \ref{sec:diff_ops_as_coinvariants}. Let us recall the construction of this action for \(N= 1\). The Lie algebroid \(\mathcal{L}_1^\sim \times_{\textnormal{Aut}_+(\cO)\ltimes K} \widehat{\fg}/\fk\) acts on the pull-back \(\mathcal{L}_1^\sim \times_{\textnormal{Aut}_+(\cO)\ltimes K} M\) of \(\Delta^\sim_K(M)\) to \(\mathcal{L}^\sim_K\). This action is compatible with the \(\widehat{\fg}_{\textnormal{out},1}\)-action. Consequently, we obtain an action of the space of differential operators on \(\mathcal{L}_1^\sim\) on \(\mathcal{L}_1^\sim \times_{\textnormal{Aut}_+(\cO)\ltimes K} M\). Since \(\mathbf{k} \in \widehat{\fg}\) acts on \(M\) by multiplication with \(k\), this action is compatible with the quantum Hamiltonian reduction, producing the desired action of \(\Delta_K(V_{K,k}) = \mathcal{D}_{K,k}\) on \(\Delta_K(M)\).

In order to describe the connection along \(T^\times C_N^X\), let us assume for simplicity that \(N = 1\) and choose a local coordinate \(u \colon U \to \bC\) for \(U \subseteq X\). Then the connection can be written as
\begin{equation}
    \nabla = \nabla_u du + \nabla_\eta d\eta,    
\end{equation}
where \(\partial_u\) is the vector field on \(X\) defined by \(u\), \(\eta\) is an exponential coordinate on the tangent direction to \(X\). Using the Segal-Sugawara operators \(S_0\) and \(S_{-1}\) defined in Section \ref{sec:segalsugawara}, the connection is now determined by
\begin{equation}
    \begin{split}
        \nabla_u = \partial_u + S_{-1} \quad \textnormal{ and }\quad \nabla_\eta = \eta\partial_\eta + S_{0};
        \end{split}
\end{equation}
see e.g.\ \cite{frenkel2004vertex,felder_kzb}.

\section{Hecke modification}\label{sec:hecke_modification}
Recall that \(X\) is a smooth complex projective curve and \(G\) is a connected simple complex algebraic group whose Lie algebra is \(\fg\). Let 
\begin{equation}
    \fg = \fn_+ \oplus \fh \oplus \fn_- = \fh \oplus \bigoplus_{\alpha \in \Phi}\fg_\alpha
\end{equation}
be a fixed triangular and root space decomposition associated to a Cartan subalgebra \(\fh \subseteq \fg\) and an associated polarized root system \(\Phi = \Phi_+ \cup \Phi_-\).

\subsection{Affine Grassmannian}
Let us start with collecting some basic facts on the affine Grassmannian, i.e.\ the ind-scheme \(\textnormal{Gr} \coloneqq \lpgrp/\lpgrpp\), following e.g.\ \cite{mirkovic_vilonen,BeilinsonDrinfeld}.

\subsubsection{Affine Grassmannian over \(X\) and \(T^\times X\).}\label{sec:affine_Grassmannian_over_X}\quad
We can put the affine Grassmannian over the curve \(X\). Indeed, \(\textnormal{Aut}(\cO)\) acts on \(\textnormal{Gr}\) by coordinate transformations, so we can define 
\begin{equation}
    \textnormal{Gr}_X \coloneqq X^\sim \times_{\textnormal{Aut}(\cO)} \textnormal{Gr} \to X.
\end{equation}
The fiber \(\textnormal{Gr}_z\) is isomorphic to \(G(\cK_z)/G(\cO_z)\), where we recall that \(\cO_z\) and \(\cK_z\) are the algebras of regular functions on the formal neighborhood \(D_z\) of \(z \in X\) and \(D_z^\circ = D_z \setminus \{z\}\) respectively.

For our purposes, it will be more reasonable to consider the larger space 
\begin{equation}
    \textnormal{Gr}_{T^\times X} \coloneqq X^\sim \times_{\textnormal{Aut}_+(\cO)} \textnormal{Gr} \to T^\times X,
\end{equation}
instead. The reason is the same as in Section \ref{sec:vbundles_hatg_modules_I}. The two spaces are related by \(\textnormal{Gr}_{T^\times X} = T^\times X \times_X \textnormal{Gr}_X\).

\subsubsection{Moduli interpretation.}\label{sec:grassmannian_as_moduli_space}\quad
The affine Grassmannian parametrizes pairs \((P,\phi)\) of \(G\)-bundles \(P\) on \(D\) with fixed trivialization \(\phi\) on \(D^\circ\). 
Similarly, the fiber \(\textnormal{Gr}_z = G(\cK_z)/G(\cO_z)\) over \(z \in X\) of the affine Grassmannian \(\textnormal{Gr}_X\) relative to \(X\) can be thought of as the moduli space of \(G\)-bundles on \(X\) with fixed trivialization over \(X^\circ_z = X \setminus\{z\}\).

\subsubsection{Orbits of the affine Grassmannian.}\label{sec:orbits_grassmannian}\quad
Let \({}^L\Lambda\) be the set of coweights of \(G\), which are the weights of the Langlands dual group \({}^LG\). Furthermore, let \({}^L\Lambda_{\ge0} \subset {}^L\Lambda\) be the subset of dominant coweights, i.e.\ those coweights \(\lambda\) satisfying \(\alpha(\lambda)\ge 0\) for all \(\alpha \in \Phi_+\). The Weyl group \(W\) of \(G\) acts on \({}^L\Lambda\) and the orbits of this action can be parametrized by \({}^L\Lambda_{\ge 0}\).
The affine Grassmannian has a decomposition 
\begin{equation}
    \textnormal{Gr} = \bigcup_{\lambda \in {}^L\Lambda/W}\textnormal{Gr}_\lambda = \bigcup_{\lambda \in {}^L\Lambda_{\ge 0}} \textnormal{Gr}_{\lambda}
\end{equation}
into \(\lpgrpp\)-orbits \(\textnormal{Gr}_\lambda \coloneqq \lpgrpp t^\lambda\) induced by 
\begin{equation}
    \lpgrp = \bigcup_{\lambda \in {}^L\Lambda_{\ge 0}}\lpgrpp t^\lambda \lpgrpp.
\end{equation}
Here, we write \(t^\lambda\) for the element of \(G(\textnormal{Spec}(\bC[t,t^{-1}])) \subseteq \lpgrp\) associated with \(\lambda \in {}^L\Lambda\), while we identify \(\lambda\) itself with an element of \(\fh\). In particular, the action of \(t^\lambda\) on \(\widehat{\fg}\) is given by
\begin{equation}\label{eq:t^lambda}
    \textnormal{Ad}_{\widehat{\fg}}(t^\lambda)a_m =\begin{cases}
         t^{\alpha(\lambda)}a_m = a_{m+\alpha(\lambda)}& ,\alpha \in \Phi,m\in\mathbb{Z},a \in \fg_\alpha
        \\ a _m + \delta_{m,0}\kappa(\lambda,a)\mathbf{k}&,m \in \mathbb{Z}, a \in \fh.
    \end{cases}
\end{equation}
This action is also referred to as ``spectral flow''.

The orbits \(\textnormal{Gr}_\lambda\) are all quasi-projective varieties and are stable under the \(\textnormal{Aut}(\cO)\)-action. Therefore, the orbit decomposition can be put over \(T^\times X\) as \(\textnormal{Gr}_{T^\times X} = \bigcup_{\lambda \in {}^L\Lambda_{\ge 0}} \textnormal{Gr}_{\lambda,T^\times X}\), where
\begin{equation}
   \textnormal{Gr}_{\lambda,T^\times X} = X \times_{\textnormal{Aut}_+(\cO)} \textnormal{Gr}_\lambda \to T^\times X.
\end{equation}
Then \(\textnormal{Gr}_{\lambda,T^\times X}\) is an algebraic variety over \(T^\times X \) whose fiber over \((z,v) \in T^\times X\) is \(G(\cO_z)s^\lambda\), where \((z,s)\) is a preimage of \((z,v)\) under \(X^\sim \to T^\times X = X^\sim/\textnormal{Aut}_+(\cO)\).

Observe that \(\textnormal{Gr}_\lambda\) can be represented as a quotient
\begin{equation}
    \textnormal{Gr}_\lambda \cong \lpgrpp/\lpgrpp_\lambda,    
\end{equation}
where \(\lpgrpp_\lambda \subseteq \lpgrpp\) is the \(\textnormal{Aut}_+(\cO)\)-invariant subgroup of \(G(\cO)\) defined by
\begin{equation}
    \lpgrpp_\lambda = \lpgrpp \cap t^{\lambda}\lpgrpp t^{-\lambda}.
\end{equation}
In particular, we observe that \(G(\cO)_\lambda\) is a nice subgroup of \(G(\cO)\) in the sense of Section \ref{sec:additional_structures} and thus we can talk about \(G(\cO)_\lambda\)-structures of \(G\)-bundles.

\subsection{Twisting \(\widehat{\fg}\)-modules} \label{sec:twisting_hatg_modules_general}
Let \(M\) be a \(\widehat{\fg}\)-module. Then for every \(g \in \lpgrp\), we can consider the \(\widehat{\fg}\)-module \(M^g\) with the same base space \(M\) but the twisted action
\begin{equation}
    a \cdot m \coloneqq (\textnormal{Ad}_{\widehat{\fg}}(g^{-1})a)m.
\end{equation}
We can think of \(M^g\) as the space of formal expressions of the form \(g m\) for \(m \in M\), because then
\begin{equation}
    a \cdot gm = g((\textnormal{Ad}_{\widehat{\fg}}(g^{-1})a)\cdot m).
\end{equation}
Therefore, we denote the identity map \(M \to M^g\) by \(m \mapsto m^g = gm \).

The modules \(M\) and \(M^g\) are isomorphic if \(g\) defines an endomorphism of \(M\) such that \(g(a\cdot m) = (\textnormal{Ad}_{\widehat{\fg}}(g)a)\cdot gm\) holds for all \(a \in \widehat{\fg}\) and \(m \in M\). For example, this is the case for \(g\in G_+(\cO)\) if all elements \(a \in \fg_+(\cO)\) act locally nilpotent on \(M\) in the sense of Section \ref{sec:exponentiation}.

Assume that \(M\) is a \((\widehat{\fg},K)\)-module for a subgroup \(K \subseteq \lpgrp\) and \(g \in \lpgrp\). Then the module \(M^g\) is integrable over \(K_g \coloneqq K\cap gKg^{-1} \subseteq K\). Furthermore, we have canonical isomorphisms \(M^g \cong M^{gk}\) and \(M^{kg} \cong M^{\textnormal{Ad}(k)g}\) for \(k \in K\).

\subsubsection{Relation to sheaves associated to \(\widehat{\fg}\)-modules.}\label{sec:twisted_Gmod_and_coinvariants}\quad
Let \(K \subseteq \lpgrpp\) be a nice subgroup, \(M\) be a \(K\)-admissible \(\widehat{\fg}\)-module of level \(k \neq -h^\vee\), and \(P\) be a \(G\)-bundle on \(X\). Choose a formal trivialization \(\sigma \in P^\sim|_{(z,s)} = \Gamma(D_z,P)\) at a point \((z,s)\in X^\sim\), where we recall that \(z \in X\) and \(s\) are formal coordinates at \(z\). The fiber of \(\Delta^{P,\sim}_K(M)\) over \([(z,s),p] \in P^\sim/(\textnormal{Aut}_+(\cO) \ltimes K)\) can be identified with \(M^g\) for the image \(g \in G(\cO)/K\) of \(p \in P^\sim|_{(z,s)}/K_z\) under the isomorphism \(P^\sim|_{(z,s)}/K_z \cong G(\cO)/K\) defined by \(\sigma\) and \(s\).

\subsubsection{Twisted vacuum module.}\label{sec:twisted_vacuum}\quad
Recall that the vacuum module 
\begin{equation}
    V_{\lpgrpp,k} = \textnormal{Ind}_{\lpalgp \oplus \bC \mathbf{k}}^{\widehat{\fg}}\bC|0\rangle    
\end{equation}
is a \((\widehat{\fg},\lpgrpp)\)-module. Up to isomorphism, the twist \(V_{\lpgrpp,k}^h\) for \(h \in \lpgrp\) depends only on the class of \(h \in \textnormal{
Gr}_\lambda\). We can collect the data of these twists in terms of a geometrization of the \emph{twisted vacuum module}
\begin{equation}
    V_{\lpgrpp,k}^\lambda \coloneqq V_{\lpgrpp,k}^{t^\lambda}.    
\end{equation}
Namely, the bundle
\begin{equation}
    \lpgrpp \times_{\lpgrpp_\lambda} V_{\lpgrpp,k}^\lambda \to \lpgrpp/\lpgrpp_\lambda \cong \textnormal{Gr}_\lambda 
\end{equation}
has the fiber \(V_{\lpgrpp,k}^h\) over \(h \in \textnormal{Gr}_\lambda\).
This can be considered as the special case of \(\Delta_{\lpgrpp_\lambda}^{P,\sim}(V_{\lpgrpp,k}^{\lambda})\), where \(P\) is the trivial bundle.

A different description of the twisted vacuum module is
\begin{equation}
    V_{G(\cO),k}^\lambda = \textnormal{Ind}_{\textnormal{Ad}(t^\lambda)\fg(\cO) \oplus \bC \mathbf{k}}^{\widehat{\fg}}\bC |0\rangle^\lambda \cong  U(t^{-1}\textnormal{Ad}_{\widehat{\fg}}(t^{\lambda})\fg[t^{-1}])|0\rangle^\lambda.
\end{equation}
In particular, for any \(a \in \fg_\alpha\) with \(\alpha \in \Phi\) and \(m \in \mathbb{Z}\) satisfying \(\alpha(\lambda) \ge m\), we have \(a_m|0\rangle^\lambda = 0\). On the other hand, for any \(h \in \fh\) and \(m \in \mathbb{Z}_{\ge 0}\), we have \(h_m|0\rangle^\lambda = -k\kappa(h,\lambda)\delta_{m,0}|0\rangle^\lambda\).

    \subsubsection{Lemma.}\label{lem:segal_sugawara_conjugation} \quad Let \(k \neq -h^\vee\) be a complex number.
    \begin{enumerate}
        \item For any \(g \in \lpgrpp\) and \(n \in \mathbb{Z}\), the Segal-Sugawara operator \(S_n\) (see Section \ref{sec:segalsugawara}) satisfies 
        \begin{equation}
            \textnormal{Ad}(g) S_n = S_n+t^{n+1}(\partial_tg)g^{-1}.
        \end{equation}

        \item For any coweight \(\lambda \in {}^L\Lambda\), we have
        \begin{equation}
            \textnormal{Ad}(t^{\lambda})S_n = S_n + \lambda_n + \delta_{n,0}\frac{k}{2}\kappa(\lambda,\lambda),
        \end{equation}
        where \(\lambda_n = t^n\lambda = t^{n+1}(\partial_t t^\lambda)t^{-\lambda}\).
        
        \item Let \(0\neq \lambda \in \fh\) be a minuscule dominant coweight, i.e.\ \(\alpha(\lambda) \in \{0,1\}\) for all positive roots \(\alpha \in \Phi_+\). Then \(V_{\lpgrpp,k}^\lambda \cong \mathbf{M}_{-\lambda^*}/N\), where \(\mathbf{M}_{-\lambda^*}\) is the \(\widehat{\fg}\)-module induced by the Verma module \(M_{-\lambda^*}\) of \(\fg\) with lowest weight \(-\lambda^* = -k\kappa(\lambda,-) \in \fh^*\) and \(N\) is the subrepresentation generated by \(e_{\theta,-1}|0\rangle^\lambda\) for the maximal root \(\theta \in \Phi_+\).
    \end{enumerate}
    \subsubsection{Proof of Lemma \ref{lem:segal_sugawara_conjugation} Part 1.}\quad
    First note that, if \emph{1}.\ holds for \(g_1,g_2 \in G(\cO)\), then it holds for \(g_1g_2 \in G(\cO)\) since
    \begin{equation}
        \begin{split}
            \textnormal{Ad}(g_1g_2)S_n &= \textnormal{Ad}(g_1)(S_n + t^{n+1}(\partial_tg_2) g_2^{-1}) \\&= S_n + t^{n+1}\textnormal{Ad}(g_1)((\partial_tg_2) g_2^{-1}) + t^{n+1}(\partial_tg_1)g_1^{-1} \\&= S_n + t^{n+1}(\partial_t (g_1g_2))(g_1g_2)^{-1}.
        \end{split} 
    \end{equation}
    Therefore, it suffices to prove \emph{1}.\ for the topological generators \(G \cup \{e^a | a\in \fg_+(\cO) = t\fg(\cO)\}\) of \(G(\cO)\).
    
    Clearly, \(\textnormal{Ad}(g)S_n = S_n\) for any \(g \in G\), since \(\textnormal{Ad}(g)\) is orthogonal with respect to the invariant bilinear form \(\kappa\). For all \(a \in \fg_+(\cO)\) we have 
    \begin{equation}
        \textnormal{Ad}(e^a)S_n = S_n + e^a[S_n,e^{-a}] = S_n + t^{n+1}\partial_ta 
    \end{equation}
    which follows from the fact that
    \begin{equation}
        [S_n,a_m] = -ma_{m+n} = -t^{n+1}\partial_ta_m.
    \end{equation}

   \subsubsection{Proof of Lemma \ref{lem:segal_sugawara_conjugation} Part 2.}
    \quad First of all, we may assume that \(\lambda\) is dominant, since \(S_n\) is invariant under the action of the Weyl group \(W \subset G\). Let us write
    \begin{equation}
        \begin{split}
        \tilde{S}_n = 2(k+h^\vee)S_n &= \sum_{i = 1}^{\textnormal{rk}(\fg)} \left(\sum_{m < 0} b_{i,m}b_{i,n-m} + \sum_{m\ge 0}b_{i,n-m}b_{i,m}\right) \\&+ \sum_{\alpha \in \Delta}\left(\sum_{m < 0} x_{\alpha,m}x_{-\alpha,n-m} + \sum_{m\ge 0}x_{-\alpha,n-m}e_{\alpha,m}\right)
        \\&\eqqcolon \tilde{S}_{n,\fh} + \tilde{S}_{n,\Delta},        
        \end{split}
    \end{equation}
    where \(\{b_i\}_{i = 1}^{\textnormal{rk}(\fg)} \subseteq \fh\) is an orthonormal basis of \(\fh\) with respect to \(\kappa\) and \(\{x_\alpha \in \fg_\alpha\}_{\alpha \in \Phi}\) are normalized such that \(\kappa(x_\alpha,x_{-\alpha}) = 1\).

    Let us assume first that \(n \neq 0\). Then using \eqref{eq:t^lambda}
    \begin{equation}
    \begin{split}
        &\textnormal{Ad}(t^{\lambda})\tilde{S}_n = \sum_{i = 1}^{\textnormal{rk}(\fg)} \left(\sum_{m < 0} b_{i,m}b_{i,n-m}+ \sum_{m\ge 0}b_{i,n-m}b_{i,m} + 2k \kappa(b_i,\lambda)b_{i,n}\right) 
        \\&+ \sum_{\alpha \in \Phi}\left(\sum_{m < 0} x_{\alpha,m+\alpha(\lambda)}x_{-\alpha,n-m-\alpha(\lambda)} + \sum_{m\ge 0}x_{-\alpha,n-m-\alpha(\lambda)}x_{\alpha,m+\alpha(\lambda)}\right) 
        \\&= \tilde{S}_{n,\fh} + 2k\lambda_n + \sum_{\alpha \in \Phi}\left(\sum_{m < \alpha(\lambda)} x_{\alpha,m}x_{-\alpha,n-m} + \sum_{m\ge \alpha(\lambda)}x_{-\alpha,n-m}x_{\alpha,m}\right)
        \\& = \tilde{S}_{n} + 2k\lambda_n + \sum_{\alpha \in \Phi_+}\sum_{m = 0}^{\alpha(\lambda)-1}[x_{\alpha,m},x_{-\alpha,n-m}] - \sum_{\alpha \in \Phi_-}\sum_{m = 1}^{-\alpha(\lambda)}[x_{\alpha,m},x_{-\alpha,n-m}]
        \\& = \tilde{S}_{n} + 2k\lambda_n + \sum_{\alpha \in \Phi}\alpha(\lambda) h_{\alpha,n} = \tilde{S}_n + 2(k+h^\vee)\lambda_n.
    \end{split}
\end{equation}
    In the last equality, we wrote \(h_\alpha \coloneqq [x_\alpha,x_{-\alpha}]\) for all \(\alpha \in \Phi\) and used that 
    \begin{equation}
        2h^\vee \lambda = \sum_{\alpha \in \Phi} [x_\alpha,[x_{-\alpha},\lambda]] = \sum_{\alpha \in \Phi} \alpha(\lambda)h_{\alpha}
    \end{equation}
    holds.

    To conclude the proof, we note that \(2S_0 = [S_1,S_{-1}]\) implies
    \begin{equation}
        \begin{split}
            \textnormal{Ad}(t^{\lambda})S_0 &= \frac{1}{2}\textnormal{Ad}(t^{\lambda})[S_1,S_{-1}] = \frac{1}{2}[S_1 + \lambda_1,S_{-1} +\lambda_{-1}] \\&= S_0 + \lambda_0 + \frac{k}{2} \kappa(\lambda,\lambda)
        \end{split}
    \end{equation}
    This concludes the proof of \emph{2.}

    \subsubsection{Proof of Lemma \ref{lem:segal_sugawara_conjugation} Part 3.}\quad  First observe that \(|0\rangle^\lambda\) is a candidate for a lowest weight vector of the \(\fg\)-action, since 
    \begin{equation}\label{eq:lowest_weight_I}
        e_{-\alpha}|0\rangle^\lambda = \textnormal{Ad}(t^{-\lambda})e_{-\alpha}|0\rangle = e_{-\alpha,\alpha(\lambda)}|0\rangle = 0    
    \end{equation}
    holds for all \(\alpha \in \Phi_+\).
    The lowest weight, that is, \(\mu \in \fh^*\) such that \(h\cdot |0\rangle^\lambda = \mu(h)|0\rangle^\lambda\), would be given by \(-\lambda ^* = -k\kappa(\lambda,-)\) because of
    \begin{equation}\label{eq:lowest_weight_II}
        h\cdot |0\rangle^{\lambda} = (\textnormal{Ad}(t^{-\lambda})h)|0\rangle = (h-\kappa(\lambda,h)\mathbf{k})|0\rangle = -\kappa(\lambda,h)k|0\rangle^\lambda.
    \end{equation}
    Now since the coweight \(\lambda\) satisfies \(\alpha(\lambda) \in \{-1,0,1\}\) for all roots \(\alpha \in \Phi\), \(\fg_+(\cO) = t\fg(\cO)\) acts trivially on the vector \(|0\rangle^\lambda\).
    Therefore, under consideration of \eqref{eq:lowest_weight_I} \& \eqref{eq:lowest_weight_II}, we have the obvious surjective morphism \(\mathbf{M}_{-\lambda^*} \to V_{G(\cO),k}^\lambda\) defined by \(|-\lambda^*\rangle \mapsto |0\rangle^\lambda\). The kernel of this morphism encodes the additional annihilation conditions \(e_{\alpha,-1}|0\rangle^\lambda = 0\) if \(\lambda(\alpha) \neq 0\) for \(\alpha \in \Phi_+\). This subrepresentation is generated by \(e_{\theta,-1}|-\lambda^*\rangle\), where \(\theta\) is the maximal root.

\subsection{Hecke stack with additional structures}
Let \(P\) be a \(G\)-bundle over \(X\) and \(z \in X\) be some point. A Hecke modification \((P',\phi)\) of \(P\) at a point \(z \in X\) consists of a \(G\)-bundle \(P'\) equipped with an isomorphism  
\begin{equation}
    \phi \colon P|_{X_z^\circ} \to P'|_{X_z^\circ}.    
\end{equation}
We now want to understand the moduli of such modifications.

\subsubsection{Idea.}\label{sec:hecke_scheme_fixed_bundle_idea}\quad
The moduli space of Hecke modifications \((P',\phi)\) at  a fixed point \(z \in X\) is given by 
\begin{equation}\label{eq:hecke_at_fixed_point}
    P^\sim|_z \times_{G(\cO_z)}G(\cK_z)/G(\cO_z) = P^\sim|_z \times_{\textnormal{Aut}(\cO) \ltimes G(\cO)} \textnormal{Gr},    
\end{equation}
where we recall that \(P^\sim|_z = \Gamma(D_z,P)\) is the set of formal trivializations of \(P\) at \(D_z\); see \cite{BeilinsonDrinfeld}. Indeed, for \(\gamma \in P^\sim|_z\) and \(\gamma' \in P^{\prime,\sim}|_z\), there exists a unique \(\widetilde{h} \in G(\cK_z)\) such that \(\phi(\gamma)\widetilde{h} = \gamma'\) and considering the \(G(\cO_z)\)-actions obtained by changing the trivializations \(\gamma\) and \(\gamma'\), we see that \(h = [\gamma,\widetilde{h}]\) defines an element of \eqref{eq:hecke_at_fixed_point}.
We can describe \(P'\) associated with this \(h = [\gamma,\widetilde{h}]\) more explicitly. The section \(\gamma\) defines a trivialization at \(D_z\) and thus identifies \(P\) with an element of \(g \in G(X_z^\circ)\!\setminus\! G(\cK_z)\). Then \(P'\) is the \(G\)-bundle defined by \(g\widetilde{h} \in G(X_z^\circ)\!\setminus\! G(\cK_z)/G(\cO_z)\). In the following, we also write \(P' = H_{z,h}(P)\).

Clearly, the notion of Hecke modification is symmetric: if \((P',\phi)\) is a Hecke modification of \(P\), then \((P,\phi^{-1})\) is a Hecke modification of \(P'\). More concretely, starting with \(h = [\gamma,\widetilde{h}] \in P^\sim|_z \times_{G(\cO_z)}\textnormal{Gr}_z\), we have that \(\gamma' = \phi(\gamma)\widetilde{h}\) defines an element of \(H_{z,h}(P)^\sim|_z\) for any preimage of \(\widetilde{h}\) in \(G(\cK)\). Then \(h^{-1} \coloneqq [\gamma',\widetilde{h}^{-1}]\) defines an element of
\(H_{z,h}(P)^\sim|_z \times_{ G(\cO)_z} \textnormal{Gr}_z\) independent of the choice of preimage of \(\widetilde{h}\). Moreover, \(H_{z,h^{-1}}(H_{z,h}(P)) = P\).

\subsubsection{Hecke scheme associated to a fixed \(G\)-bundle.}\label{sec:hecke_scheme_fixed_bundle}\quad
We have seen in Section \ref{sec:orbits_grassmannian} that \(\textnormal{Gr} = \bigcup_{\lambda \in {}^L\Lambda_{\ge0}} \textnormal{Gr}_\lambda\) and it was noted in Section \ref{sec:vbundles_hatg_modules_I} that, in general, we should only quotient out the action of \(\textnormal{Aut}_+(\cO) \subset\textnormal{Aut}(\cO)\). Therefore, we can consider Hecke modifications of type \(\lambda\), that is, those which are parameterized by
\begin{equation}
    P^\sim|_z \times_{\textnormal{Aut}_+(\cO) \ltimes \lpgrpp} \textnormal{Gr}_{\lambda}
\end{equation}
and vary them over \(z \in X\) and extend it to \(T^\times X\) in order to obtain:
\begin{equation}
    \begin{split}
        \textnormal{Hecke}^P_{\lambda} &\coloneqq P^\sim \times_{\textnormal{Aut}_+(\cO) \ltimes \lpgrpp} \textnormal{Gr}_{\lambda} \\&\cong P^\sim/(\textnormal{Aut}_+(\cO) \ltimes \lpgrpp_\lambda) = C^P_{\lpgrpp_\lambda}.
    \end{split}
\end{equation}
Here, we used that \(\textnormal{Gr}_\lambda \cong \lpgrpp/\lpgrpp_\lambda\). 

More generally, we may allow \(\underline{K}\)-structures, for a family \(\underline{K} = (K_1,\dots,K_N)\) of nice subgroups of \(G(\cO)\), at points different from the point of Hecke insertion. In this way, we recover the moduli space 
\begin{equation}
    \textnormal{Hecke}^P_{\underline{K},\lambda} \coloneqq C^P_{(\underline{K},\lpgrpp_\lambda)}.
\end{equation}

Let us mention how the symmetry of the Hecke modification carries over. In Section \ref{sec:hecke_scheme_fixed_bundle_idea}, we have seen that to any \(h \in P^\sim|_z \times_{G(\cO)} \textnormal{Gr}_z\) we can associate an \(h^{-1} \in H_{z,h}(P)^\sim|_z \times_{G(\cO)} \textnormal{Gr}_z\). From the construction via inversion inside \(G(\cK)\), it is clear that for \(((z,v),h) \in \textnormal{Hecke}^P_{\underline{K},\lambda}\) this construction defines an element 
\begin{equation}\label{eq:inverse_of_h}
    ((z,v),h^{-1}) \in \textnormal{Hecke}^{H_{z,h}(P)}_{\underline{K},-\lambda} = \textnormal{Hecke}^{H_{z,h}(P)}_{\underline{K},\lambda^\circ}.
\end{equation}
Here, \(\lambda^\circ = -w_0(\lambda)\) is the dual coweight of \(\lambda\), i.e.\ the unique dominant coweight in the Weyl group orbit of \(-\lambda\), whereby \(w_0\) is the longest element of the Weyl group \(W\) of \(G\).

\subsubsection{General Hecke 
stack.}\label{sec:general_Hecke}\quad
    The Hecke scheme can be generalized in order to also let the \(G\)-bundle \(P\) vary. Indeed, we can adapt the Hecke stack from \cite{BeilinsonDrinfeld}, to the setting of \(G\)-bundles with additional structures. 
    Throughout the remainder of this section, \(\underline{K} = (K_1,\dots,K_N)\) is a family of nice subgroups of \(\lpgrpp\).
    
    The Hecke stack of coweight \(\lambda \in {}^L\Lambda_{\ge0}\) with additional \(\underline{K}\)-structures at \(N\) marked points is the substack \(\textnormal{Hecke}^X_{\underline{K},\lambda} \subseteq \textnormal{Bun}_{G,\underline{K}}^X \times_{T^\times C^X_{N}} \textnormal{Bun}_{G,(\underline{K},\lpgrpp)}^X\) consisting of pairs \((P_1,P_2)\) of \(G\)-bundles with \(\underline{K}\)-structures at \(N\)-marked points, which are isomorphic outside a \((N+1)\)-st marked point of \(X\), with respect to an isomorphism of type \(\lambda\).
    Fixing a \(G\)-bundle \(P\) on \(X\) with \(\underline{K}\)-structures at the first \(N\) marked points, we recover \(\textnormal{Hecke}^P_{\underline{K},\lambda}\).

    A convenient way to describe \(\textnormal{Hecke}^X_{\underline{K},\lambda}\) is given as follows:
    \begin{equation}\label{eq:hecke_general}
        \begin{split}
            \textnormal{Hecke}^X_{\underline{K},\lambda} &= \textnormal{Bun}_{G,(\underline{K},\{e\})}^{X} \times_{T^\times C_{N+1}^X,\lpgrpp} (\textnormal{Gr}_{\lambda,T^\times X} \times_{T^\times X} T^\times C^X_{N+1}) \\&\cong \textnormal{Bun}^{X}_{G,(\underline{K},\lpgrpp_\lambda)}.
        \end{split}
    \end{equation}
    Here, we used the fact that \(\textnormal{Gr}_\lambda \cong \lpgrpp/\lpgrpp_\lambda\) holds and used the projection \(\pi_{N+1}\colon T^\times C^X_{N+1} \to T^\times X\) given by 
    \(((z_1,v_1),\dots,(z_{N+1},v_{N+1})) \mapsto (z_{N+1},v_{N+1})\) in order to define the fiber product \(\textnormal{Gr}_{\lambda,T^\times X} \times_{T^\times X} T^\times C^X_{N+1}\).

    Observe that the Hecke stack carries a natural involution that exchanges \(P_1\) and \(P_2\). This involution restricts to an identification 
    \begin{equation}\label{eq:sigma_lambda}
        \sigma_\lambda \colon \textnormal{Hecke}^X_{\underline{K},\lambda} \to \textnormal{Hecke}^X_{\underline{K},\lambda^\circ}.    
    \end{equation}
    In the identification \eqref{eq:hecke_general}, one can express this involution using the inversion of Hecke parameters. Let us write \([P,((z,v),h)] \in \textnormal{Hecke}_{\underline{K},\lambda}^X\) with \(P  \in \textnormal{Bun}^X_{G,\underline{K}}\) and \(((z,v),h) \in \textnormal{Hecke}_\lambda^P\). Then using \eqref{eq:inverse_of_h}, we obtain
    \begin{equation}
        [P,((z,v),h)] \stackrel{\sigma_\lambda}\mapsto [H_{z,h}(P),((z,v),h^{-1})].
    \end{equation}

   \subsection{Hecke modification}
    We can now define a map 
    \begin{equation}\label{eq:hecke_on_Gbundles}
        \eta_\lambda\colon \textnormal{Hecke}_{\underline{K},\lambda}^X  \to \textnormal{Bun}_{G,\underline{K}}^X\,,\qquad [P,((z,v),h)] \mapsto H_{z,h}(P). 
    \end{equation}
    Here,  \(P \in \textnormal{Bun}^X_{G,\underline{K}}\), \([(z,v),h] \in \textnormal{Hecke}^P_\lambda\) for a \(z \in X\) different from the marked points of \(X\) at which \(P\) has additional structures, and \(H_{z,h}(P)\) was defined in Section \ref{sec:hecke_scheme_fixed_bundle_idea}.
    
    \subsubsection{Hecke modification of coinvariants I: Motivation.}\label{sec:hecke_fixed_bundle}\quad 
    In \cite{BeilinsonDrinfeld}, the Hecke modification for \(\lambda \in {}^L\Lambda_{\ge 0}\) is defined as a \(\mathcal{D}\)-module pull-push along the correspondence
    \begin{equation}
        \textnormal{Bun}_{G}^X \leftarrow \textnormal{Hecke}_\lambda^X \to X \times \textnormal{Bun}_{G}^X\,,\qquad P \mapsfrom (P,P',z,\phi) \mapsto (z,P').
    \end{equation}
    We might want to adapt this by defining the Hecke modification of \(\Delta_{\underline{K}}(\underline{M})\) as the pull-push of twisted \(\mathcal{D}\)-modules via the analogous correspondence 
    \begin{equation}
        \textnormal{Bun}_{G,\underline{K}}^X \leftarrow \textnormal{Hecke}_{\underline{K},\lambda}^X\rightarrow \textnormal{Bun}_{G,(\underline{K},\lpgrpp)}^X.    
    \end{equation}
    After applying the natural involution \((P,P',z,\phi) \mapsto (P',P,z,\phi^{-1})\) to both maps, which does not change the pull-push functor, we may identify the map along which we pull-back with the Hecke map \(\eta_\lambda\) from \eqref{eq:hecke_on_Gbundles}.
    Since outside of the critical level there are no global twisted differential operators, we only do the pull-back along \(\textnormal{Bun}_{G,\underline{K}}^X \stackrel{\eta_\lambda}\leftarrow \textnormal{Hecke}_{\underline{K},\lambda}^X\) and do not want to consider the push-forward afterwards, i.e.\ we keep the additional structure at the new marked point. 

    \subsubsection{Hecke modification of coinvariants II: Definition.}\label{sec:hecke_def}\quad
    Let \(\underline{M} = (M_1,\dots,M_N)\) be a family of \(\underline{K}\)-admissible \(\widehat{\fg}\)-modules of level \(k \neq -h^\vee\).
    Following the intuition outlined in Section \ref{sec:hecke_fixed_bundle}, we define the Hecke modification on coinvariants by
    \begin{equation}
        H_\lambda(\Delta_{\underline{K}}(\underline{M})) = \eta^*_\lambda \Delta_{\underline{K}}(\underline{M})
    \end{equation}
    on the level of quasi-coherent sheaves. By dualizing, this defines the Hecke modification for conformal blocks as well. 
    
    In particular, this means that 
    \begin{equation}
        H_\lambda(\Delta_{\underline{K}}(\underline{M}))\big|_{C^P_{\underline{K},((z,v),h)}} = \Delta_{\underline{K}}^{H_{z,h}(P)}(\underline{M})\big|_{C^{P,\circ}_{\underline{K},z}}
    \end{equation}
    holds, consistent with \cite{jeong_lee_nekrasov,teschner_notes}. Here, \(C^{P,\circ}_{\underline{K},z} \subset C_{\underline{K}}^P\) is the preimage under \(C^P_{\underline{K}} \to C^X_N\) of the open subset of points \((z_1,\dots,z_N) \in C^X_N\) satisfying \(z_i \neq z\) for all \(1 \le i \le N\) and \(C^{P,\circ}_{\underline{K},z}\cong C^P_{\underline{K},((z,v),h)} \subset \textnormal{Hecke}_{\underline{K},\lambda}^X\) is the subset where the underlying \(G\)-bundle \(P\) and the Hecke parameter \(((z,v),h) \in \textnormal{Hecke}_\lambda^P\) are fixed.

\section{Main results}\label{sec:results}
As before, \(X\) is a smooth complex projective curve and \(G\) is a connected simple complex algebraic group whose Lie algebra is \(\fg\). Moreover, we fix a non-critical level \(k \in \mathbb{C}\), \(k \neq -h^\vee\).

\subsection{Realizing Hecke modifications}
The first set of results of this paper clarifies the twisted \(\mathcal{D}\)-module structure of Hecke modified coinvariants and their construction using insertions of twisted vacuum modules. More precisely, we have the following theorem.

    \subsubsection{Theorem.}\quad \label{thm:Hecke_and_insertion}
    Let \(\lambda\) be a dominant coweight of \(G\), \(\underline{K} = (K_1,\dots,K_N)\) be a family of nice subgroups of \(G(\cO)\), and \(\underline{M}\) be a \(\underline{K}\)-admissible family of \(\widehat{\fg}\)-modules of level \(k\).

    \begin{enumerate}
        \item There is a canonical isomorphism 
        \begin{equation}\label{eq:Hecke_is_insertion}
            H_\lambda(\Delta_{\underline{K}}(\underline{M})) \stackrel{\cong}\to  \Delta_{(\underline{K},G(\cO)_\lambda)}(\underline{M},V_{G(\cO),k}^\lambda)
        \end{equation}
        of quasi-coherent sheaves.
        
        \item There is a natural \(\eta_\lambda\)-morphism between the twisted differential operators on \(\textnormal{Hecke}_{\underline{K},\lambda}^X\) and \(\B_{G,\underline{K}}^X\) in the sense of \cite[Section 1.4.]{beilinson_bernstein_jantzen_conjecture}, so the pull-back \(H_\lambda(\Delta_{\underline{K}}(\underline{M})) =  \eta_\lambda^*\Delta_{\underline{K}}(\underline{M})\) has a natural twisted \(\mathcal{D}\)-module structure over \(\textnormal{Hecke}_{\underline{K},\lambda}^X\) by virtue of \cite[Section 1.4.]{beilinson_bernstein_jantzen_conjecture}. Furthermore, \eqref{eq:Hecke_is_insertion} becomes an isomorphism of twisted \(\mathcal{D}\)-modules over \(\textnormal{Hecke}_{\underline{K},\lambda}^X = \B_{G,(\underline{K},G(\cO)_\lambda)}^X\).
    \end{enumerate}

    \subsubsection{Proof of Theorem \ref{thm:Hecke_and_insertion} Part 1.}\quad
    In order to understand what happens at the \((N+1)\)-th point, i.e.\ the point of Hecke modification, we first recall that inserting a vacuum module does not change the space of coinvariants. In particular, we have 
    \begin{equation}
        \Pi^{*}\Delta_{\underline{K}}(\underline{M}) \cong \Delta_{(\underline{K},\lpgrpp_{\lambda^\circ})}(\underline{M},V_{G(\cO),k}),
    \end{equation}
    where \(\Pi \colon \textnormal{Hecke}_{\underline{K},\lambda^\circ}^X = \B_{G,(\underline{K},G(\cO)_{\lambda^\circ})}^X \to \B_{G,\underline{K}}^X\) is the canonical projection, forgetting about the \((N+1)\)-th marked point and the additional structure at that point. Here, recall that \(\lambda^\circ = -w_0\lambda\) is the dual coweight of \(\lambda\), where \(w_0\) is the longest element of the Weyl group of \(G\). Now we have \(\eta_\lambda = \Pi\sigma_\lambda\) for 
    \begin{equation}
        \sigma_\lambda\colon \textnormal{Hecke}_{\underline{K},\lambda}^X \to
        \textnormal{Hecke}_{\underline{K},\lambda^\circ}^X\,,\qquad [P,((z,v),h)] \mapsto [H_{z,h}(P),((z,v),h^{-1})] 
    \end{equation}
    defined in Section \ref{sec:general_Hecke}.
    Therefore, we have
    \begin{align}\label{eq:Hecke_via_sigma}
        &H_\lambda(\Delta_{\underline{K}}(\underline{M})) = \eta_\lambda^*\Delta_{\underline{K}}(\underline{M}) = \sigma_{\lambda}^*\Delta_{(\underline{K},\lpgrpp_{\lambda^\circ})}(\underline{M},V_{G(\cO),k}).
    \end{align}
    It remains to show that
    \begin{equation}\label{eq:Hecke_is_insertion}
            \sigma_{\lambda}^*\Delta_{(\underline{K},\lpgrpp_{\lambda^\circ})}(\underline{M},V_{G(\cO),k}) \cong  \Delta_{(\underline{K},G(\cO)_\lambda)}(\underline{M},V_{G(\cO),k}^\lambda)
        \end{equation}
    holds.

    The sheaves \(\Delta_{(\underline{K},G(\cO)_{\lambda^\circ})}(\underline{M},V_{G(\cO),k})\) and \(\Delta_{(\underline{K},G(\cO)_\lambda)}(\underline{M},V_{G(\cO),k}^\lambda)\) descend from trivial bundles on the scheme
    \begin{equation}
        \textnormal{Hecke}_{\mu,N}^{X,\sim} = \B_{G,N+1}^{X,\sim} \times G(\cO)t^\lambda G(\cO)
    \end{equation}
    by modding out the \(\prod_{i = 1}^{N+1}(\textnormal{Aut}_+(\cO) \ltimes K_i)\) and \(\widehat{\fg}_{\textnormal{out},N+1}\)-actions.
    Here, \(\mu \in \{\lambda,\lambda^\circ\}\), \(K_{N+1} = G(\cO)\), and this group acts by conjugation on \(G(\cO)t^\lambda G(\cO)\) while all other \(\textnormal{Aut}_+(\cO) \ltimes K_i\) act trivially on this orbit.
    Furthermore, \(\sigma_\lambda\) lifts to a morphism 
     \begin{equation}
        \sigma_\lambda^\sim\colon \textnormal{Hecke}_{\underline{K},\lambda}^{X,\sim} \to
        \textnormal{Hecke}_{\underline{K},\lambda^\circ}^{X,\sim}
    \end{equation}
    which is defined on the fiber over \(((z_1,s_1),\dots,(z_{N+1},s_{N+1})) \in C_{N+1}^{X,\sim}\) as
    \begin{equation}
        ([g_1,\dots,g_{N+1}],h) \mapsto ([g_1,\dots,g_{N+1}h],h^{-1}),
    \end{equation}
    where \([g_1,\dots,g_{N+1}] \in G(X_{\underline{z}}^\circ)\!\setminus\!\prod_{i = 1}^{N+1}G(\cK_{\underline{z}})\) and \(h \in G(\cO)t^{\lambda}G(\cO)\).

    It remains to define a \(\sigma_\lambda^\sim\)-equivariant isomorphism  
    \begin{equation}
        \textnormal{Hecke}_{\underline{K},\lambda^\circ}^{X,\sim} \times (M_1\otimes \dots \otimes M_N \otimes V_{G(\cO),k}) \to
        \textnormal{Hecke}_{\underline{K},\lambda}^{X,\sim} \times (M_1\otimes \dots \otimes M_N \otimes V_{G(\cO),k}^\lambda)
    \end{equation}
    that is compatible with the respective \(\prod_{i = 1}^{N+1}(\textnormal{Aut}_+(\cO) \ltimes K_i)\) and \(\widehat{\fg}_{\textnormal{out},N+1}\)-actions.
    The desired map is given on the fiber over \(((z_1,s_1),\dots,(z_{N+1},s_{N+1})) \in C_{N+1}^{X,\sim}\) as 
    \begin{equation}\label{eq:twisting_of_vacuum_lift}
        \begin{split}  
            (([g_1,\dots,g_{N+1}],h^{-1}=g_-^{-1}t^{-\lambda}g_+^{-1}),m \otimes g_+^{-1}v) \mapsto (([g_1,\dots,g_{N+1}h^{-1}],h),m \otimes t^{\lambda}g_-v),
        \end{split} 
       \end{equation}
       where \(m \in M_1 \otimes \dots \otimes M_N\), \(g_1,\dots,g_{N+1} \in G(\cO)\), \(g_+,g_- \in G(\cO)\), and \(v \in V_{G(\cO),k}\).
    It is straight-forward to check that this is an \(\sigma_\lambda^\sim\)-equivariant map compatible with the aforementioned actions, concluding the proof.\qed

    \subsubsection{Proof of Theorem \ref{thm:Hecke_and_insertion} Part 2.}\quad
    Recall that
    \begin{equation}
        \begin{split}
            &\mathcal{D}_{\underline{K},k} = \Delta_{\underline{K}}(V_{K_1,k},\dots,V_{K_N,k})
        \end{split}
    \end{equation}
    holds. Therefore, \eqref{eq:Hecke_is_insertion} provides the identification
    \begin{equation}
        \eta^*_\lambda\mathcal{D}_{\underline{K},k} = \Delta_{(\underline{K},G(\cO)_\lambda)}(V_{K_1,k},\dots,V_{K_N,k},V_{G(\cO),k}^\lambda).    
    \end{equation}
    This space of coinvariants has a natural action of \(\mathcal{D}_{(\underline{K},G(\cO)_\lambda),k}\)-action compatible with its right \(\eta^{-1}_\lambda\mathcal{D}_{\underline{K},k}\)-action. Therefore, \cite[Section 1.4.]{beilinson_bernstein_jantzen_conjecture} provides \(H_\lambda(\Delta_{\underline{K}}(\underline{M})) = \eta_\lambda^*\Delta_{\underline{K}}(\underline{M})\) with a twisted \(\mathcal{D}\)-module structure over \(\textnormal{Hecke}_{\underline{K},\lambda}^X\). The \(\mathcal{D}_{(\underline{K},G(\cO)_\lambda),k}\)-action is thereby defined as the natural action on
    \begin{equation}
        \eta^*_\lambda\Delta_{\underline{K}}(\underline{M}) = \cO_{\textnormal{Hecke}_{\underline{K},\lambda}^X} \otimes_{\eta^{-1}\cO_{\B_{G,\underline{K}}^X}} \Delta_{\underline{K}}(\underline{M}) \cong \eta^*_\lambda\mathcal{D}_{\underline{K},k} \otimes_{\eta_\lambda^{-1}\mathcal{D}_{\underline{K},k}}\Delta_{\underline{K}}(\underline{M}).
    \end{equation}
    The definition of the \(\mathcal{D}\)-module structure along \(C_{N+1}^X\) uses a similar strategy.

    Clearly, the identification \eqref{eq:Hecke_is_insertion} is a \(\mathcal{D}_{(\underline{K},G(\cO)_\lambda),k}\)-module morphism. 
    The connection over \(C^X_{N+1}\) on the image of sections of \(\eta_\lambda^*\Delta_{\underline{K}}(\underline{M})\) under the isomorphism \eqref{eq:Hecke_is_insertion} acts by \(S_0\) and \(S_{-1}\) conjugated by \(h\) at precisely the \((N+1)\)-th marked point; see \eqref{eq:twisting_of_vacuum}. This coincides with the action of those elements on sections of \(\Delta_{(\underline{K},\lpgrpp_\lambda)}(\underline{M},V_{G(\cO),k}^\lambda)\), which concludes the proof of Part \emph{2.}\ of the theorem. \qed

    \subsubsection{Remark: description of \eqref{eq:Hecke_is_insertion} for fixed \(G\)-bundles.}\quad 
    Let us describe the isomorphism \eqref{eq:Hecke_is_insertion} explicitly fiberwise for any fixed \(G\)-bundle \(P\). We have the identification
    \begin{equation}\label{eq:twisting_of_vacuum}
        \begin{split}   \Delta_{\lpgrpp_{\lambda^\circ}}^{H_{z,h}(P),\sim}(V_{G(\cO),k})\Big|_{[(z,v),h^{-1}]} &\stackrel{\cong}\to \Delta_{\lpgrpp_\lambda}^{P,\sim}(V_{G(\cO),k}^\lambda)\Big|_{[(z,v),h]}
        \\
        [((z,u),p),m] &\mapsto [((z,u),p\widetilde{h}^{-1}),\widetilde{h}m].
        \end{split} 
       \end{equation}
    Here, \((z,u)\in X^\sim\) is a preimage of \((z,v) \in T^\times X \cong X^\sim/\textnormal{Aut}_+(\cO)\), \(p \in H_{z,h}(P)^\sim|_z\) is a preimage of \(h^{-1} \in H_{z,h}(P)^\sim|_z/G(\cO_z)_{\lambda^\circ}\), \(m \in V_{G(\cO),k}\), and \(\widetilde{h} \in G(\cO) t^\lambda G(\cO)\) is a preimage of the \(\textnormal{Gr}_{\lambda,z}\)-component of \(h\). The map \eqref{eq:twisting_of_vacuum} is well-defined, since \(V_{G(\cO),k}\) is integrable over \(G(\cO)\) and recalling the geometrization of twisted vacuum modules from Section \ref{sec:twisted_vacuum}. 

     On the other hand,
    \begin{equation}\label{eq:proof_hecke_is_insertion_I}
        \Delta^{H_{z,h}(P),\sim}_{\underline{K}}(\underline{M})\big|_{C_{\underline{K},z}^{P,\circ}} \cong \Delta_{\underline{K}}^{P,\sim}(\underline{M})\big|_{C_{\underline{K},z}^{P,\circ}}
    \end{equation}
    holds, since \(H_{z,h}(P)\) is isomorphic to \(P\) outside of \(z\). Here, \(C^P_{\underline{K},((z,v),h)} \subset \textnormal{Hecke}_{\underline{K},\lambda}^X\) and \(C_{\underline{K},z}^{P,\circ} \subset C_{\underline{K}}^P\) were defined at the end of Section \ref{sec:hecke_fixed_bundle}. 
    
    Combining \eqref{eq:twisting_of_vacuum} and \eqref{eq:proof_hecke_is_insertion_I}, we obtain
    \begin{equation}\label{eq:proof_hecke_is_insertion_II}
        \Delta^{H_{z,h}(P),\sim}_{(\underline{K},G(\cO)_{\lambda^\circ})}(\underline{M},V_{G(\cO),k})\big|_{C_{\underline{K},z}^{P,\circ}} \cong \Delta_{(\underline{K},G(\cO)_\lambda)}^{P,\sim}(\underline{M},V_{G(\cO),k}^\lambda)\big|_{C_{\underline{K},z}^{P,\circ}}.
    \end{equation}
    The identification \eqref{eq:proof_hecke_is_insertion_II} is equivariant with respect to the embedding 
    \begin{equation}
        \widehat\fg_{\textnormal{out},N}^P \to \widehat\fg_{\textnormal{out},N+1}^{H_{z,h}(P)}|_{z_{N+1} = z}    
    \end{equation}
    where \((-)|_{z_{N+1} = z}\) means that the \((N+1)\)-th point of \(X\) is fixed to \(z\).

\subsection{Description via coordinate transformations}\label{sec:birational}
It is possible to describe the Hecke modification of coinvariants with \(N\) inserted modules, which are coinvariants with \((N+1)\) inserted modules according to Theorem \ref{thm:Hecke_and_insertion} via appropriate coordinate transformation. This was used for \(\fg = \mathfrak{sl}_2(\mathbb{C})\) and \(X= \mathbb{P}^1\) in \cite{jeong_lee_nekrasov,teschner_notes}, where this transformation becomes very explicit, as we will see in Section \ref{sec:explicit_calculations}. 
We consider this phenomenon in two parts, first for \(X = \mathbb{P}^1\), where for a fixed Hecke parameter, we can even control the isomorphism class of \(G\)-bundles under the Hecke transformation, and then for a general higher genus curve \(X\).

\subsubsection{The \(X = \mathbb{P}^1\) case I.}\quad
Let us choose a global coordinate \(t\) of \(\mathbb{A}^1 = \mathbb{P}^1\setminus\{\infty\}\). Recall that any \(G\)-bundle \(P\) on \(\mathbb{P}^1\) is given, up to isomorphism, by some coweight \(\mu\) by gluing the trivial \(G\)-bundle on \(\mathbb{A}^1\) and \(\mathbb{A}^1_{\infty} = \mathbb{P}^1\setminus\{0\}\) along \(\mathbb{A}^1\setminus\{0\}\) using \(t^\mu \in G(\textnormal{Spec}(\bC[t,t^{-1}]))\). This amounts to the description 
\begin{equation}
    \begin{split}
        \B_G^{\mathbb{P}^1} = G(\mathbb{A}^1_\infty) \!\setminus\! G(\mathbb{A}^1\setminus\{0\})/G(\mathbb{A}^1)  \cong {}^L\Lambda /W \cong {}^L\Lambda_{\ge 0}
    \end{split}
\end{equation}
of the moduli space of \(G\)-bundles, where \(W\) is the Weyl group of \(G\). More formally, the points \(\mu \in {}^L\Lambda_{\ge 0}\) are understood as the stacky points \(\{P_\mu\}/\textnormal{Aut}(P_\mu)\),
where \(P_\mu\) is the \(G\)-bundle on \(\mathbb{P}^1\) determined by \(\mu \in {}^L\Lambda_{\ge 0}\).

For any \(\lambda \in {}^L\Lambda_{\ge 0}\), we also have
\begin{equation}
    \textnormal{Gr}_\lambda \cong \lpgrpp/\lpgrpp_\lambda \cong G(\mathbb{A}^1)/G(\mathbb{A}^1)_\lambda,
\end{equation}
where \(G(\mathbb{A}^1)_\lambda = G(\mathbb{A}^1) \cap t^{\lambda}G(\mathbb{A}^1)t^{-\lambda}\).
Consequently, we can write
\begin{equation}\label{eq:section_of_orbit}
     T^\times \mathbb{A}^1 \times \textnormal{Gr}_\lambda \stackrel{\cong}\to \textnormal{Gr}_{\lambda,\mathbb{P}^1}|_{T^\times \mathbb{A}^1} \,,\qquad ((z,v),[h(t)]) \mapsto ((z,v),[h(t-z)]).
\end{equation}
An open subset \(U\subseteq \textnormal{Gr}_{\lambda,\mathbb{P}^1}\) can now be described by means of the open subset \(N_+^\lambda t^\lambda \subseteq\textnormal{Gr}_\lambda\), where using
\begin{equation}
     \fn_+^\lambda \coloneqq \lpalgp \cap \textnormal{Ad}(t^{\lambda})t^{-1}\fg[t^{-1}] \subseteq \fn_+(\cO) = \fn_+ \otimes \cO    
\end{equation}
we define \(N_+^\lambda \coloneqq \exp(\fn_+^\lambda) \subseteq \lpgrpp\).

\subsubsection{Proposition.}\label{prop:hecke_modifications_isomorphism_class}\quad
Let \(\lambda\) and \(\mu\) be dominant coweights of \(G\). The Hecke modification of the \(G\)-bundle defined by \(\mu\) with respect to 
\begin{equation}
    ((z,v),\exp(xt^j))\in T^\times \mathbb{A}^1 \times N_+^\lambda t^\lambda \cong \textnormal{Gr}_{\lambda,\mathbb{P}^1}|_{\mathbb{A}^1},
\end{equation}
for \(\alpha \in \Phi_+\) satisfying \(0\le j<\alpha(\lambda)\) and \(x \in \fg_\alpha\), is isomorphic to the \(G\)-bundle defined by the coweight \(\nu= \mu + \lambda -(j+\alpha(\mu))\alpha^\vee\). In particular, the isomorphism class does not depend on \((z,v) \in T^\times \mathbb{A}^1\).

\subsubsection{Proof of Proposition \ref{prop:hecke_modifications_isomorphism_class}.}\quad\label{sec:proof_of_prop_2}
For \(t_z = t-z\), the transition map of the Hecke modified \(G\)-bundle is
\begin{equation}
    t_z^\mu \exp(xt_z^j)t_z^\lambda = \exp(xt_z^{j+\alpha(\mu)})t_z^{\mu + \lambda} = t_z^{\mu + \lambda}\exp(xt_z^{j-\alpha(\lambda)}).
\end{equation}
We need to factorize this element into the form \(A(t_z)t_z^{\eta}B(t_z)\) with \(A\) being regular at \(t_z = \infty\) and \(B\) being regular at \(t_z = 0\) in order to determine the isomorphism class. 

Let us start with the \(\textnormal{rk}(G) = 1\) case, i.e.\ $G=\textnormal{SL}_2(\mathbb{C})$ or \(G = \textnormal{PGL}_2(\mathbb{C})\). We want 
to factorize
\begin{equation}\label{eq:rk1_factorization}
\begin{split}
    &\begin{pmatrix}
      1 & at_z^{j+2p}\\
      0 & 1
  \end{pmatrix}\begin{pmatrix}
      t_z^{q+p} & 0 \\ 0 & t_z^{-q-p}
  \end{pmatrix} = \begin{pmatrix}
      t_z^{q+p} & at_z^{j+p-q}\\0 & t_z^{-p-q}
  \end{pmatrix} \\&= \begin{pmatrix}
      t_z^{q+p} & 0 \\ 0 & t_z^{-q-p}
  \end{pmatrix}  \begin{pmatrix}
      1 & at_z^{j-2q}\\
      0 & 1
  \end{pmatrix}
  = A(t_z) \begin{pmatrix}
      t_z^{r} & 0 \\ 0 & t_z^{-r}
  \end{pmatrix}B(t_z)
\end{split}
\end{equation}
for every \(a \in \bC\), where $A(t_z)$ is regular at
$t_z = \infty$ and $B(t_z)$ is regular at $t_z = 0$.
Here, \(2p = \alpha(\mu),2q = \alpha(\lambda),\) and \(2r = \alpha(\eta)\) for the unique positive root \(\alpha\). If \(G = \textnormal{SL}_2(\mathbb{C})\), then \(p,q,r \in \mathbb{Z}\), but if \(G = \textnormal{PGL}_2(\mathbb{C})\), they can be elements of \(\frac{1}{2}\mathbb{Z}\).

We have obvious solutions of \eqref{eq:rk1_factorization} for 
$A$ or $B$ except if $-2p < j < 2q$.  In this case, we use the identity
\begin{equation}
    \begin{split}
        &\begin{pmatrix}
      t_z^{q+p} & at_z^{j+p-q}\\0 & t_z^{-p-q}
  \end{pmatrix} = \begin{pmatrix}
      0 & t_z^{-q+p+j}\\-t_z^{q-p-j} & a^{-1}t_z^{-p-q}
  \end{pmatrix}\begin{pmatrix}
       a^{-1}& 0\\t^{-j+2q}_z & a
  \end{pmatrix}
  \\&=\begin{pmatrix}
      0 & 1\\ -1& a^{-1}t_z^{-j-2p}
  \end{pmatrix}\begin{pmatrix}
      t_z^{q-p-j} & 0 \\ 0 & t_z^{-q+p+j}
  \end{pmatrix}
  \begin{pmatrix}
      a^{-1} & 0\\t_z^{-j+2q} & a
  \end{pmatrix}.
    \end{split}
\end{equation}
With the Weyl and Cartan group elements
\begin{equation}
    w = \begin{pmatrix}
    0 & 1 \\ -1 &0
\end{pmatrix}\quad \textnormal{ and }\quad h_a = \begin{pmatrix}
    a & 0\\ 0 & a^{-1}
\end{pmatrix}
\end{equation}
respectively, this provides
\begin{equation}
    A(t_z)=w   \begin{pmatrix}
      1 & -a^{-1}t_z^{-j-2p}\\ 0 & 1
  \end{pmatrix}\quad\textnormal{ and }\quad B(z)= \begin{pmatrix}
      1 &  0\\ at_z^{-j+2q} & 1
  \end{pmatrix}h_a^{-1},
\end{equation}
as solutions for the factorization problem \eqref{eq:rk1_factorization} in the case \(\textnormal{rk}(G) = 1\).

In the case \(\textnormal{rk}(G) > 1\), we want to factorize $\exp(xt_z^{j+\alpha(\mu)})t_z^{\mu+\lambda}$ with $x\in\mathfrak{g}_\alpha$. Consider the subalgebra of \(\fg\) isomorphic to \(\mathfrak{sl}_2(\mathbb{C})\) spanned by $e_{\pm\alpha}\in \mathfrak g_{\pm\alpha}$
and $\alpha^\vee\in\mathfrak h$ satisfying $[e_\alpha,e_{-\alpha}]=\alpha^\vee$ and $[\alpha^\vee,e_{\pm\alpha}]=\pm
2e_{\pm\alpha}$. The corresponding subgroup $G_\alpha\in G$ admits a canonical surjective homomorphism
$\textnormal{SL}_2(\mathbb{C}) \to G_\alpha$. Let $w_\alpha$ be the image of $w$ under this homomorphism. Then we claim that the factorization for $x = ae_\alpha$
is solved by
\begin{equation}
    \exp(ae_\alpha t_z^{j+\alpha(\mu)})t_z^{\mu+\lambda}=A(t_z)t_z^{\nu}B(t_z),
\end{equation}
where \( \nu = \mu + \lambda - (j+\alpha(\mu))\alpha^\vee \) and
    \begin{equation}
             A(t_z)=w_\alpha \exp(-a^{-1}e_\alpha t_z^{-j-\alpha(\mu)}) \quad \textnormal{ and } \quad 
             B(t_z)=\exp(ae_{-\alpha}t_z^{\alpha(\lambda)-j})a^{-\alpha^\vee}
        \end{equation}
for the Weyl reflection \(w_\alpha\) defined by \(\alpha\).

To prove this we set $p=\alpha(\mu)/2$ and $q = \alpha(\lambda)/2$. These are integers if $G_\alpha=\textnormal{SL}_2(\mathbb{C})$ and may be half-integers if
$G_\alpha=\textnormal{PGL}_2(\mathbb{C})$.  Notice that 
we have 
\begin{equation}
    \begin{split}
        &[(p+q)\alpha^\vee-\mu -\lambda,e_{-\alpha}] = -2(p+q)+\alpha(\mu)+\alpha(\lambda)  = 0 
        \\&\implies [(p+q)\alpha^\vee - \mu -\lambda,B(t_z)] = 0.
    \end{split}
\end{equation}
Therefore, the claimed factorization is equivalent to
\begin{equation}
    \exp(ae_\alpha t_z^{j+\alpha(\mu)})t_z^{(p+q)\alpha^\vee}=A(t_z)t_z^{(p-j-q)\alpha^\vee}B(t_z),
\end{equation}
which takes place in $G_\alpha$ and reduces to the \(\textnormal{rk}(G) = 1\) calculation above. \qed

\subsubsection{Remark to Proposition \ref{prop:hecke_modifications_isomorphism_class}.}\quad
In general, note that the coweight \(\lambda + \mu - (j + \alpha(\mu))\alpha^\vee\) does not need to be dominant. Therefore, the \(G\)-bundle determined by \(\lambda + \mu - (j + \alpha(\mu))\alpha^\vee\) is actually the \(G\)-bundle determined by the unique dominant coweight in the Weyl group orbit of this coweight.

\subsubsection{The \(X = \mathbb{P}^1\) case II.}\quad
Let \(\lambda\) and \(\mu\) be dominant coweights and \((\alpha,j) \in \Phi_+ \times \mathbb{N}_0\) be a pair such that \(j < \alpha(\lambda)\). 
Consider the \(G\)-bundles \(P\) and \(Q\) on \(\mathbb{P}^1\) determined by \(\mu\) and 
\begin{equation}
    \nu = \mu + \lambda - (j+\alpha(\mu))\alpha^\vee    
\end{equation}
respectively. 
According to Proposition \ref{prop:hecke_modifications_isomorphism_class}, we have an isomorphism of the Hecke modification of \(P\) with respect to the Hecke parameter \(((z,v),u_{\alpha,j} \coloneqq \exp(u t^j)) \in T^\times\mathbb{A}^1 \times N_+^\lambda t^\lambda\), where \(u \in \fg_\alpha\), and \(Q\). This isomorphism is given by a change of trivialization on \(D_z\) and by left multiplication with \(h_{\alpha,j}(t-z)^{-1}\) for 
\begin{equation}\label{eq:gauge}
       h_{\alpha,j}(t) =
           w_\alpha \exp(-u^{-1}e_\alpha t^{-j-\alpha(\mu)})t^{\lambda}.
\end{equation}
We obtain an induced map
\begin{equation}\label{eq:definition_p_ualphaj}
    \begin{split}
        &p_{u_{\alpha,j}} \colon  C_{\underline{K},u_{\alpha,j}}^P \coloneqq \bigsqcup_{(z,v) \in T^\times\mathbb{A}^1}C_{\underline{K},((z,v),u_{\alpha,j})}^{P} \to C_{\underline{K},\infty}^{Q,\circ}
    \end{split}
\end{equation}
using the notations from the end of Section \ref{sec:hecke_fixed_bundle}. More precisely, \(p_{u_{\alpha,j}}\) is defined on the fiber of \(C_{\underline{K},((z,v),u_{\alpha,j})}^{P}\) over \(((z_1,v_1),\dots,(z_N,v_N),(z,v)) \in C^{\mathbb{A}^1}_{N+1}\) via
\begin{equation}\label{eq:def_p}
    \begin{split}
    c \coloneqq [g_1,\dots,g_N,u_{\alpha,j}] \mapsto [h_{\alpha,j}(t_1-z)^{-1}g_1,\dots,h_{\alpha,j}(t_N-z)^{-1}g_N] \eqqcolon p_{u_{\alpha,j}}(c).
    \end{split}
\end{equation}
The fiber of \(\Delta_{\underline{K}}^{Q,\sim}(\underline{M})\) over \(p(c)\) is spanned by tensors
\begin{equation}\label{eq:image_of_Phi_0}
    T = m_1\otimes \dots \otimes m_N \in M_{1,z_1}^{h_{\alpha,j}^{-1}g_1} \otimes \dots \otimes M_{N,z_N}^{h_{\alpha,j}^{-1}g_N}.    
\end{equation}
By assigning to this tensor the element
\begin{equation}\label{eq:image_of_Phi}
    h_{\alpha,j}(\bullet - z)T \otimes |0\rangle^h\coloneqq h_{\alpha,j}(t_1 - z)m_1 \otimes \dots \otimes h_{\alpha,j}(t_N - z)m_N \otimes |0\rangle^{h_{\alpha,j}},
\end{equation}
which is in the fiber \(
    M^{g_1}_{1,z_1} \otimes \dots \otimes M^{g_N}_{N,z_N} \otimes V_{\lpgrpp,k,z}^{h_{\alpha,j}}    \)
of
\(\Delta_{(\underline{K},\lpgrpp_\lambda)}^{P,\sim}(\underline{M},V_{\lpgrpp,k}^\lambda)\) over \(c\), we have defined a map
\begin{equation}\label{eq:Phi_definition_precoinvariants}
        p_{u_{\alpha,j}}^*\Delta^{Q,\sim}_{\underline{K}}(\underline{M}) \to \Delta_{(\underline{K},\lpgrpp_\lambda)}^{P,\sim}(\underline{M},V_{\lpgrpp,k}^\lambda)\Big|_{C_{(\underline{K},\lpgrpp_\lambda)}^{P,u_{\alpha,j}}}. 
    \end{equation}
Let us remark that \eqref{eq:image_of_Phi} is defined using that
\begin{equation}
    h_{\alpha,j}(t_i - z) \in G(\mathbb{A}^1\setminus\{z\}) \subseteq G(\cO_{z_i})
\end{equation}
holds for \(z_i \notin \{z,\infty\}\), so \(h_{\alpha,j}(t_i-z)\) indeed defines an isomorphism \(M^{h^{-1}_{\alpha,j}g_i}_{i,z_i} \to M^{g_i}_{i,z_i}\) of \((\widehat{\fg}_{z_i},K_{i,z_i})\)-modules.

\subsubsection{Theorem.}\label{thm:N_to_N+1_P^1case}\quad 
    The map \eqref{eq:Phi_definition_precoinvariants} induces a unique isomorphism 
    \begin{equation}\label{eq:Phi_full_def}
        \Phi\colon p_{u_{\alpha,j}}^*\Delta^Q_{\underline{K}}(\underline{M}) \to H_\lambda(\Delta^P_{\underline{K}}(\underline{M}))\Big|_{C_{\underline{K},u_{\alpha,j}}^{P}} 
    \end{equation}
    of twisted \(\mathcal{D}\)-modules.

\subsubsection{Proof of Theorem \ref{thm:N_to_N+1_P^1case}.}\quad\label{proof:thm_N_to_N+1_P^1case}
Let us first proof that \eqref{eq:Phi_full_def} defines an isomorphism of vector bundles.
Let \(\eta_{\alpha,j}\) be the restriction of \(\eta_\lambda\) from \eqref{eq:hecke_on_Gbundles} to \(C_{\underline{K},u_{\alpha,j}}^P\). Then 
\begin{equation}
    H_\lambda(\Delta^P_{\underline{K}}(\underline{M}))\Big|_{C_{\underline{K},u_{\alpha,j}}^{P}}  = \eta_{\alpha,j}^*\Delta_{\underline{K}}^P(\underline{M}).
\end{equation}
Clearly, \(p_{\alpha,j} = \iota_{\alpha,j}\eta_{\alpha,j}\), where 
\begin{equation}
    \iota_{\alpha,j}\colon \bigsqcup_{(z,v)\in T^\times \mathbb{A}^1}C_{\underline{K},z}^{H_{z,u_{\alpha,j}}(P),\circ} \to C_{\underline{K},\infty}^{Q,\circ}
\end{equation}
is the map defined on the fiber of  \(C_{\underline{K},z}^{H_{z,u_{\alpha,j}}(P),\circ}\) over \(((z_1,v_1),\dots,(z_N,v_N)) \in C^{\mathbb{A}^1}_{N+1}\) via
\begin{equation}\label{eq:def_iota}
    \begin{split}
    [g_1,\dots,g_N] 
    \mapsto [h_{\alpha,j}(t_1-z)^{-1}g_1,\dots,h_{\alpha,j}(t_N-z)^{-1}g_N].
    \end{split}
\end{equation}
It remains to prove that the assignment
\begin{equation}\label{eq:Phi_definition_precoinvariants_2}
    m_1\otimes \dots \otimes m_N \mapsto h_{\alpha,j}(t_1 - z)m_1 \otimes \dots \otimes h_{\alpha,j}(t_N - z)m_N,
\end{equation}
where \(m_i \in M_i^{h_{\alpha,j}^{-1}g_i}\),
defines a bundle isomorphism
\begin{equation}\label{eq:iota_and_coinvariants}
    \iota_{\alpha,j}^*\Delta^{Q}_{\underline{K}}(\underline{M})\big|_{C_{\underline{K},z}^{H_{z,u_{\alpha,j}}(P),\circ}} \cong\Delta_{\underline{K}}^{H_{z,u_{\alpha,j}}(P)}(\underline{M})\big|_{C_{\underline{K},z}^{H_{z,u_{\alpha,j}}(P),\circ}}
\end{equation}
for every \((z,v) \in T^\times \mathbb{A}^1\). But this immediately follows from the fact that \eqref{eq:Phi_definition_precoinvariants_2} is equivariant under
\begin{equation}
    \textnormal{Ad}(h_{\alpha,j}(\bullet - z)) \colon \widehat{\fg}_{\textnormal{out},N}^Q \to \widehat{\fg}_{\textnormal{out},N}^{H_{z,u_{\alpha,j}}(P)},
\end{equation}
since 
\begin{equation}\label{eq:equivariance_under_gout}
    {\small\begin{split}
        &\textnormal{Ad}(h_{\alpha,j}(\bullet-z)a)\cdot (h_{\alpha,j}(t_1 - z)m_1 \otimes \dots \otimes h_{\alpha,j}(t_N - z)m_N) 
        \\ &= \sum_{i = 1}^N h_{\alpha,j}(t_1-z)m_1 \otimes \dots \otimes h_{\alpha,j}(t_i -z)am_i\otimes \dots \otimes h_{\alpha,j}(t_N-z)m_N
    \end{split}}
\end{equation}
holds.

Next, we investigate the compatibility of the twisted \(\mathcal{D}\)-module structure. By virtue of \eqref{eq:iota_and_coinvariants} and the statements in Section \ref{sec:diff_ops_as_coinvariants}, we have
\begin{equation}
    \iota_{\alpha,j}^*\mathcal{D}_{\underline{K},k}^Q\big|_{C_{\underline{K},z}^{H_{z,u_{\alpha,j}}(P),\circ}} \cong \mathcal{D}_{\underline{K},k}^{H_{z,u_{\alpha,j}}(P)}\big|_{C_{\underline{K},z}^{H_{z,u_{\alpha,j}}(P),\circ}}.    
\end{equation}
Consequently, \cite{beilinson_bernstein_jantzen_conjecture} equips \(p_{\alpha,j}^*\Delta_{\underline{K}}^Q(\underline{M})\) with a natural twisted \(\mathcal{D}\)-module structure along the additional structures of \(Q\). Furthermore, this twisted \(\mathcal{D}\)-module structure is now clearly compatible with the bundle isomorphism \eqref{eq:Phi_full_def}, using a similar calculation as \eqref{eq:equivariance_under_gout}.
 
Finally, we consider the connection along \(C_{N+1}^X\). Using the notation from \eqref{eq:image_of_Phi_0} \& \eqref{eq:image_of_Phi}, we have 
\begin{equation}
    \begin{split}
        &(\partial_{t_i} + S_{-1}^{(i)})(h_{\alpha,j}(\bullet - z)T \otimes |0\rangle^h)\\&= 
        h_{\alpha,j}(\bullet -z)((\partial_{t_i} + S_{-1}^{(i)})T) \otimes |0\rangle^h+ \textnormal{Ad}(h_{\alpha,j}(t_i - z))(\partial_{t_i} + S_{-1}^{(i)})T \otimes |0\rangle^h
        \\& =h_{\alpha,j}(\bullet -z)((\partial_{t_i} + S_{-1}^{(i)})T) \otimes |0\rangle^h.
    \end{split} 
\end{equation}
for \(i \in \{1,\dots,N\}\), where we used Lemma \ref{lem:segal_sugawara_conjugation}. Furthermore,
\begin{equation}
    \begin{split}
        &(\partial_{z} + S_{-1}^{(N+1)}) (h_{\alpha,j}(\bullet - z)T \otimes |0\rangle^h)  \\&= h_{\alpha,j}(\bullet-z)\partial_{z}T \otimes |0\rangle^{h_{\alpha,j}} + 
        \textnormal{Ad}(h_{\alpha,j}(t_i - z))((\partial_{z} + S_{-1,N+1})^{(N+1)})(T \otimes |0\rangle^{h_{\alpha,j}})
        \\& = \left(h_{\alpha,j}(\bullet-z)(\partial_z + (\partial_zh_{\alpha,j}(\bullet-z))h_{\alpha,j}(\bullet-z)^{-1} T\right) \otimes |0\rangle^{h_{\alpha,j}},
    \end{split} 
\end{equation}
where in the bracket in the last expression, we indeed find \(h_{\alpha,j}(\bullet-z)\) applied to the derivation of \(T\) with respect to the connection in the \((N+1)\)-th coordinate induced by the pull-back via \(p_{\alpha,j}\). \qed

\subsubsection{Remark.}\quad As we will see explicitly in Section \ref{sec:explicit_calculations}, the multiplication by a rational map \(X \to G\) like \(h_{\alpha,j}\) that is used in \(\Phi\) can be made explicit in the form of coordinate transformations if we can identify the modules \(M_i\) with \(\mathcal{D}\)-modules on appropriate flag varieties, as is the case for Verma modules. This is the reason we have referred to Theorem \ref{thm:N_to_N+1_P^1case} as a description through coordinate transformations before. 

\subsubsection{Example: minuscule weight and trivial bundle.}\label{exp:N_to_N+1_P^1case}\quad
In general, we cannot simply provide a description of coinvariants as in Theorem \ref{thm:N_to_N+1_P^1case} while varying the Hecke parameter over \(U\), since this changes the isomorphism class of \(Q\). 
A special case where this is not an issue is given by \(\mu = 0\), i.e.\ \(P = P_{\textnormal{triv}}\) is the trivial \(G\)-bundle, and miniscule \(\lambda\). Indeed, then \(j = 0 = \alpha(\mu)\) and so \(\nu = \lambda\) in Proposition \ref{prop:hecke_modifications_isomorphism_class}. In particular, then \(Q\) is the \(G\)-bundle on \(\mathbb{P}^1\) defined by \(\lambda\).

Therefore, we vary the parameter \(u_{\alpha,j}\) in \eqref{eq:definition_p_ualphaj} in oder to obtain a map
\begin{equation}
    \begin{split}
        &p_{u_{\alpha,j}}\colon  C_{\underline{K},U}^{P} \to C_{\underline{K},\infty}^{Q,\circ},
    \end{split}
\end{equation}
where \(C_{\underline{K},U}^{P}\subseteq C_{(\underline{K},\lpgrpp_\lambda)}^{P}\) consists of those elements 
\begin{equation}
    [((z_1,v_1),g_1),\dots,((z_{N+1},v_{N+1}),g_{N+1})] \in C_{(\underline{K},\lpgrpp_\lambda)}^{P} 
\end{equation}
for which \(z_1,\dots,z_N,z = z_{N+1} \in \mathbb{A}^1\) and \(g_{N+1} \in U\).

    Similarly, the map \(\Phi\) from Theorem \ref{thm:N_to_N+1_P^1case} can now be varied over the \(u_{\alpha,j}
    \in U\) to obtain a map
    \begin{equation}
        \Phi\colon p_{u_{\alpha,j}}^*\Delta^Q_{\underline{K}}(\underline{M}) \to H_\lambda(\Delta^P_{\underline{K}}(\underline{M}))\Big|_{C_{\underline{K},U}^{P}} 
    \end{equation}
    of coinvariants compatible with the canonical connections.

\subsubsection{The general case.}\qquad
Recall that any point of \(X\) admits an affine open neighborhood \(Y = X \setminus\{y_1,\dots,y_m\}\) with a local coordinate \(s\), i.e.\ an \'etale map \(s\colon Y \to \bC\).  
Furthermore, let us choose a section \(h\) of \(G(\mathbb{A}^1)t^\lambda G(\mathbb{A}^1) \to \textnormal{Gr}_\lambda\), e.g.\ an element of \(U = N_+^\lambda t^\lambda\). 
Then
\begin{equation}
    h_z \coloneqq h(s-s(z)) \colon Y \to G\,,\qquad w\mapsto h(s(w)-s(z))
\end{equation}
defines an element of \(G(Y\setminus\{z\})\). 

Fix a \(G\)-bundle \(P\) that is trivialized over \(Y\) and let \(Q \coloneqq H_{z,h_z}(P)\). The \(G\)-bundle \(Q\) is obtained from \(P\) by adding \(h_z\) as a transition map at \(z\). We obtain an induced map
\begin{equation}\label{eq:definition_p_ualphaj}
    \begin{split}
        &p_{((z,v),h_z)} \colon  C_{\underline{K},((z,v),h_z)}^P \to C_{\underline{K},z}^{Q,\circ}
    \end{split}
\end{equation}
defined on the fiber of \(C_{\underline{K},((z,v),h_z)}^{P}\) over \(((z_1,v_1),\dots,(z_N,v_N),(z,v)) \in C^Y_{N+1}\) via
\begin{equation}\label{eq:def_p}
    \begin{split}
    c \coloneqq [g_1,\dots,g_N,h_z] \mapsto [h_z^{-1}g_1,\dots,h_{z}^{-1}g_N] \eqqcolon p_{((z,v),h_z)}(c).
    \end{split}
\end{equation}
The fiber of \(\Delta_{\underline{K}}^{Q,\sim}(\underline{M})\) over \(p_{((z,v),h_z)}(c)\) is spanned by tensors
\begin{equation}\label{eq:image_of_Phi_0}
    T = m_1\otimes \dots \otimes m_N \in M_{1,z_1}^{h_{z}^{-1}g_1} \otimes \dots \otimes M_{N,z_N}^{h_{z}^{-1}g_N}.    
\end{equation}
By assigning to this tensor the element
\begin{equation}\label{eq:image_of_Phi}
    h_{z}(T \otimes |0\rangle)\coloneqq h_{z}m_1 \otimes \dots \otimes h_{z}m_N \otimes |0\rangle^{h_{z}},
\end{equation}
which is in the fiber \(
    M^{g_1}_{1,z_1} \otimes \dots \otimes M^{g_N}_{N,z_N} \otimes V_{\lpgrpp,k,z}^{h_z}    \)
of
\(\Delta_{(\underline{K},\lpgrpp_\lambda)}^{P,\sim}(\underline{M},V_{\lpgrpp,k}^\lambda)\) over \(c\), we have defined
\begin{equation}\label{eq:Phi_definition_precoinvariants}
        p_{((z,v),h_z)}^*\Delta^{Q,\sim}_{\underline{K}}(\underline{M}) \to \Delta_{(\underline{K},\lpgrpp_\lambda)}^{P,\sim}(\underline{M},V_{\lpgrpp,k}^\lambda)\Big|_{C_{\underline{K},((z,v),h_z)}^{P}}. 
    \end{equation} 
Using now essentially the same arguments as in the Proof \ref{proof:thm_N_to_N+1_P^1case}, we can see that the following holds.
\subsubsection{Theorem.}\label{thm:N_to_N+1_general_case}\quad 
    The map \eqref{eq:Phi_definition_precoinvariants} induces a unique isomorphism 
    \begin{equation}
        \Phi\colon p_{((z,v),h_z)}^*\Delta^{Q}_{\underline{K}}(\underline{M}) \to H_\lambda(\Delta_{\underline{K}}^{P}(\underline{M}))\Big|_{C_{\underline{K},((z,v),h_z)}^{P}}. 
    \end{equation}
    of coinvariants compatible with the canonical connections.

\subsubsection{Outlook.}\quad It would be interesting to see if the isomorphism in Theorem \ref{thm:N_to_N+1_general_case} can be varied along the Hecke parameter.
For instance, this would be possible, at least locally, if \(H_{z,h_z}(P)\) are isomorphic to each other for \((z,h_z)\) in an open subset of \(\textnormal{Gr}_{\lambda,T^\times Y}\). \'Etale locally, this is for instance possible for all \(G\)-bundles satisfying  \(H^1(X,\fg^P) = 0\).
Another open question would be if the \(G\)-bundles \(P\) can be varied, i.e.\ if \(\Phi\) could be generalized to a substack of \(\B_{G,\underline{K}}^X\).

\section{Explicit calculations for \(G = \textnormal{PGL}_2(\mathbb{C})\)} \label{sec:explicit_calculations}

\subsection{Generalities on parabolic structures}

\subsubsection{Conformal blocks for parabolic structures.}\quad
Any \(G\)-module \(P\) on \(X\) admits a reduction to the Borel subgroup \(B\subseteq G\); see \cite{drinfeld_simpson}. Therefore, if we denote by \(\mathbf{B} \subseteq \lpgrpp\) the subgroup consisting of \(g \in \lpgrpp\) with \(g(0) \in B\), we have that \(C^P_{(\mathbf{B},\dots,\mathbf{B})} \cong C^X_N  \times (G/B)^N\). Now the coinvariants and conformal blocks associated to a family of \((\widehat{\fg},\mathbf{B})\)-modules \(\underline{\mathbf{M}} = (\mathbf{M}_1,\dots,\mathbf{M}_N)\) on \(C^P_{(\mathbf{B},\dots,\mathbf{B})}\) induced by \((\fg,B)\)-modules \(\underline{M} = (M_1,\dots,M_N)\) can be, fiberwise over \(C_N^X\), described using the \(\mathcal{D}\)-modules on \(G/B\) associated with these modules in the Beilinson-Bernstein theory. This will be used in the following calculations.

\subsubsection{\(\mathcal{D}\)-module description of contragradient Verma modules.}\quad 
We want to understand the Beilinson-Bernstein localization of contragradient Verma modules explicitly. Therefore,
let us identify the large cell in \(N_- \subseteq G/B\) with \(\mathbb{A}^{|\Delta_-|}\) and chose root coordinates on \(\{y_\alpha\}_{\alpha \in \Delta_-}\) of \(N_-\).

According to e.g.\ \cite{frenkel_langlads_correspondance_for_koop_groups}, the contragredient Verma module of \(\fg\) with highest weight \(\chi \in \fh^*\) can now be identified with the space of functions \(M_\lambda^* = \bC[N_-] = \bC[(y_\alpha)_{\alpha \in \Delta_-}]\) equipped with maps \(\rho_\chi \colon \fg \to \textnormal{Der}(\bC[N_ -])\) which are defined on standard generators \(\{e_i,h_i,f_i\}_{i = 1}^{\textnormal{dim}(\fh)}\) via
\begin{equation}
    \begin{split}
        &\rho_\chi(e_i) = \partial_{\alpha_i} + \sum_{\alpha \in \Delta_-} P^{(i)}_\alpha \partial_\alpha,\\
        &\rho_\chi(h_i) = \sum_{\alpha \in \Delta_-} \alpha(h_i) y_\alpha \partial_\alpha + \chi(h_i),\\
        &\rho_\chi(f_i) = \sum_{\alpha \in \Delta_-}Q^{(i)}_\alpha\partial_\alpha + \chi(h_i)y_\alpha.
    \end{split} 
\end{equation}
Here, \(P^{(i)}_\alpha,Q^{(i)}_\beta \in \bC[N_-]\) are some polynomials and we wrote \(\partial_\alpha = \partial/\partial y_\alpha\). 

\subsubsection{The case \(X = \mathbb{P}^1\) and \(G = \textnormal{PGL}_2(\bC)\).}\quad For \(G = \textnormal{PGL}_2(\bC)\) the above steps can be made more explicit: 
\begin{enumerate}
    \item We have
    \begin{equation}
        G/B = \left\{\begin{bmatrix}
        * & a \\ * & b
    \end{bmatrix}\Bigg| [a:b] \in \mathbb{P}^1 \right\} \cong \mathbb{P}^1
    \end{equation}
    and \(N_- \cong \mathbb{A}^1 = \mathbb{P}^1\setminus \{\infty\}\) is obtained by \(a \neq 0\).
    \item The \(\mathcal{D}\)-modules to  \(M_\chi\) on \(G/B\) coincide with \(\chi\)-differentials \(\Omega^\chi\). After trivializing over \(N_- = \mathbb{A}^1\), the differential operators associated with \(\{e,h,f\} \subset \fg = \mathfrak{sl}_2(\bC)\) on \(M_\chi^* = \bC[\mathbb{A}^1] = \bC[y]\) are given by
    \begin{equation}
        \begin{split}
            &\rho_\chi(e) = \partial\\
            &\rho_\chi(h) = -2y \partial + \chi\\
            &\rho_\chi(f) = -y^2\partial + \chi y,
        \end{split} 
    \end{equation}
    where \(\chi = \chi(h) \in \bC\).

    \item \(G = \textnormal{PSL}_2(\bC)\) acts on \(G/B\) by Möbius transformations. Indeed,  
    \begin{equation}
        \begin{pmatrix}
            m_1 & m_2 
            \\ m_3 &m_4  
        \end{pmatrix}
        \begin{pmatrix}
            * & a 
            \\ * & b
        \end{pmatrix} = \begin{pmatrix}
            * & m_1a+m_2b \\
            * & m_3a+m_4b
        \end{pmatrix}.
    \end{equation}
    Consider for example the matrices \(\begin{pmatrix}
        -x & 1 \\t-z& 0
        \end{pmatrix}\) parametrized by \(x \in \mathbb{A}^1\). 
    Then the associated coordinate transformations are given by
    \begin{equation}\label{eq:coordinate_transform}
        y = \begin{pmatrix}
            1 \\ y
        \end{pmatrix} \mapsto \eta \coloneqq \begin{pmatrix}
        -x & 1 \\t-z& 0
        \end{pmatrix}\begin{pmatrix}
            1 \\ y
        \end{pmatrix} = \begin{pmatrix}
           y-x \\ t-z
        \end{pmatrix} = -\frac{t-z}{x-y}.
    \end{equation}
    Therefore, the associated bundle has the gluing map \(y \mapsto \eta = -\frac{t-z}{x-y}\) and we have that in \(\Omega^{\chi}\) 
    \begin{equation}
        f(y) \mapsto \left(\frac{(x-y)^2}{t-z}\right)^\chi f\left(\eta\right)
    \end{equation}
    holds.
\end{enumerate}

\subsection{Calculations}
\subsubsection{Setup for explicit calculation.}\label{sec:explicit_calcs_setup}\quad Let us now turn to the explicit calculations. A section \(\Psi\) of the conformal blocks \(\mathcal{C}^{P_{\textnormal{triv}}}_{(\mathbf{B},\dots,\mathbf{B})}(\mathbf{M}_{\chi_1},\dots,\mathbf{M}_{\chi_N})\) for the trivial \(\textnormal{PGL}_2(\mathbb{C})\)-bundle \(P_{\textnormal{triv}}\) over \((\mathbb{A}^N \times N_-^N) \cap C^{P_{\textnormal{triv}}}_{(\mathbf{B},\dots,\mathbf{B})} \subset C^{P_{\textnormal{triv}},\circ}_{(\mathbf{B},\dots,\mathbf{B})}\) is a polynomial in 2N-variables
\begin{equation}
    \Psi \in \bC[x_1,\dots,x_N;t_1,\dots,t_N]
\end{equation}
satisfying the \(\Gamma(\mathbb{P}^1,\fg_{\textnormal{out}}) = \fg = \mathfrak{sl}_2(\bC)\) invariance
\begin{equation}
    \begin{aligned}
    0 &= \rho_\chi(h)\Psi = \rho_\chi(e)\Psi = \rho_\chi(f)\Psi
    \\&=\sum_{i=1}^N(-2x_i \partial_{x_i}+2\chi_i)\Psi=  \sum_{i=1}^N\partial_{x_i}\Psi=
    \sum_{i=1}^N(-x_i^2\partial_{x_i}+2\chi_ix_i)\Psi.  
    \end{aligned}
\end{equation}
These equations are also called Ward identities.

Moreover, \(\Psi\) satisfies the KZ-equation, which state that it defines a flat section with respect to the KZ-connection. More precisely,
\begin{equation}
    (k+2)\frac{\partial}{\partial t_i}\Psi(x_1,...,x_N;t_1,...,t_N)=\sum^N_{\substack{j=1\\ i \neq j}} \frac{\Omega_{ij}}{t_i-t_j}\Psi(x_1,...,x_N;t_1,...,t_N)
\end{equation}
 where $\Omega_{ij} = e^{(i)}f^{(j)} + f^{(i)}e^{(j)} + \frac{1}{2}h^{(i)}h^{(j)} \in U(\mathfrak{sl}_2(\bC))^{\otimes N}$.

\subsubsection{Goal of calculation.} \quad We want to verify the abstract statement of Theorem \ref{thm:N_to_N+1_P^1case} explicitly. Since we work with flag varieties, it is however more convenient to work with antidominant coweights, which simply amounts to a Weyl group conjugation. We consider the antidominant minuscule weight
\begin{equation}
    \lambda = \frac{1}{2}\begin{pmatrix}
        -1&0\\0&1
    \end{pmatrix} \implies t^\lambda = \begin{pmatrix}
    t^{-1/2} & 0\\ 0 & t^{1/2} 
\end{pmatrix}  = \begin{pmatrix}
    1 & 0\\ 0 & t 
\end{pmatrix} \in \textnormal{PGL}_2(\bC(\!(t)\!))
\end{equation}
and the trivial \(G\)-bundle \(P = P_{\textnormal{triv}}\), so we are (up to exchanging \(\lambda\) with \(-\lambda\)) in the setting of Example \ref{exp:N_to_N+1_P^1case}.
The orbit \(\textnormal{Gr}_\lambda\) is isomorphic to \(G/B\) and so we may take \(U = N_- \cong \mathbb{A}^1\) as the open subset parametrized by \(\begin{pmatrix}
        1 & 0 \\y^{-1}& 1
        \end{pmatrix}\) and \eqref{eq:gauge} used in Theorem \ref{thm:N_to_N+1_P^1case} is given by
\begin{equation}
    \begin{split}
    h = \begin{pmatrix}
        -y & 1\\t & 0
        \end{pmatrix}.
    \end{split}
\end{equation}
The action \(h(\bullet-t_{N+1})\) on the \(i\)-th factor is now precisely given by the change of variables $x_i \mapsto \xi_i=-\frac{t_i-t_{N+1}}{x_i-x_{N+1}}$ as we have seen in \eqref{eq:coordinate_transform}. Therefore, we want to verify that 
\begin{equation}
    \Upsilon=\prod_{i=1}^N \left(  \frac{(x_i-x_{N+1})^2}{(t_i-t_{N+1})} \right)^{\chi_i} \Psi(\xi_1,...,\xi_N;t_1,...,t_N)
\end{equation}
defines a section of 
\begin{equation}\label{eq:isomorphism_parabolic_structures}
    H_\lambda \mathcal{C}^Q_{(\mathbf{B},\dots,\mathbf{B})}(\mathbf{M}_{\chi_1},\dots,\mathbf{M}_{\chi_N}) \cong \mathcal{C}^Q_{(\mathbf{B},\dots,\mathbf{B},\mathbf{B})}(\mathbf{M}_{\chi_1},\dots,\mathbf{M}_{\chi_N},\mathbf{M}_{k/2}).  
\end{equation}
Here, we used the fact that the twisted vacuum module is isomorphic to the Verma module of highest weight \(k/2\) (see Lemma \ref{lem:segal_sugawara_conjugation} and change the sign of \(\lambda\)) and \(Q\) corresponds to the rank two bundle \(\mathcal{O}_{\mathbb{P}^1} \oplus \mathcal{O}_{\mathbb{P}^1}(1)\). 

Let us make the claim that \(\Upsilon\) is a section of \eqref{eq:isomorphism_parabolic_structures} more concrete. Our invariance algebra is given by the global sections of \(\widehat{\fg}^Q\), where \(Q\) is the \(\textnormal{PGL}_2(\bC)\)-bundle associated to \(\mathcal{O}_{\mathbb{P}^1} \oplus \mathcal{O}_{\mathbb{P}^1}(1)\). Therefore, this invariance algebra is given by 
\begin{equation}
    \Gamma(\mathbb{P}^1,\widehat{\fg}_{\textnormal{out}}^Q) = \bC (h - \mathbf{k})\oplus \bC e  \oplus t\bC e,
\end{equation}
where the shift by the central element in the Cartan component is obtained through \(
    t^{-\lambda}\partial_t t^{\lambda } =  \lambda/t
    \)
and \(\textnormal{Ad}_{\widehat{\fg}}(t^{\lambda})(xt^m) = \textnormal{Ad}_{\lpalg}(xt^{m}) + \textnormal{res}_0\kappa(t^{-\lambda}\partial_t t^{\lambda },xt^m)\mathbf{k}\).

Therefore, section \(\Upsilon\) of \eqref{eq:isomorphism_parabolic_structures} by definition have to satisfy
\begin{equation}
    \begin{aligned}
    0 = \sum_{i=1}^{N+1}(-2x_i \partial_{x_i}+2\chi_i)\Upsilon - k\Upsilon= \sum_{i=1}^{N+1}\partial_{x_i} \Upsilon = \sum_{i = 1}^{N+1}t_i\partial_{x_i}\Upsilon, 
    \end{aligned}
\end{equation}
where \(\chi_{N+1} = k/2\).

Moreover, we want to explicitly recognize that indeed the $(N+1)$-point KZ-equations 
\begin{equation}
    (k+2)\frac{\partial}{\partial t_i}\Upsilon=\sum^{N+1}_{\substack{j=1\\ i \neq j}} \frac{\Omega_{ij}}{t_i-t_j}\Upsilon 
\end{equation}
for $i=1,...,N$ are satisfied. 

\subsubsection{Explicit calculation of the Ward identities.}\quad
In the new coordinates, the Ward identities for \(\Psi\) take the form

\begin{equation}\label{eq:Ward_psi_in_xi}
    \begin{aligned}
    0 = \sum_{i=1}^N(-2\xi_i \partial_{\xi_i}+2\chi_i)\Psi=
    \sum_{i=1}^N\partial_{\xi_i}\Psi=
    \sum_{i=1}^N(-\xi_i^2\partial_{\xi_i}+2\chi_i\xi_i)\Psi.  
    \end{aligned}
\end{equation}
Let us use these to derive the Ward identities for $\Upsilon$. We write $\Upsilon=A\Psi$ with \(A = \prod_{i = 1}^{N}\frac{(x_i-x_{N+1})^{2\chi_i}}{(t_i-t_{N+1})^{\chi_i}}\) and immediately observe that
\begin{equation}\label{eq:first_Ward_for_Y}
    \partial_{x_{N+1}}\Upsilon=  -\sum_{i=1}^{N}\partial_{x_i}\Upsilon \iff \sum_{i = 1}^{N+1}\partial_{x_i}\Upsilon = 0
\end{equation}
holds, which is the first Ward identity for \(\Upsilon\).

Let us observe that for any test function \(\psi = \psi(\xi_1,\dots,\xi_N;t_1,\dots,t_N)\) and writing 
\begin{equation}
    \eta(x_1,\dots,x_{N+1};t_1,\dots,t_{N+1}) = A\psi    
\end{equation}
we have the following identities
\begin{equation}\label{eq:product_rule}
    \begin{split}
    &\partial_{x_i} \eta = A\partial_{x_i}\psi + (\partial_{x_i}A) \eta =  \frac{(t_i-t_{N+1})}{(x_i-x_{N+1})^2}A \partial_{\xi_i}\psi+\frac{2\chi_i}{x_i-x_{N+1}}\eta\\
    &\implies \partial_{\xi_i} \psi= \frac{(x_i-x_{N+1})^2}{(t_i-t_{N+1})}A^{-1}\left(\partial_{x_i}\eta-\frac{2\chi_i}{x_i-x_{N+1}}\eta\right).
    \end{split}
\end{equation}
This implies in particular 
\begin{equation}\label{eq:lie_action_coordinate_change}
    \begin{split}
    &(-2\xi_i\partial_{\xi_i}+2\chi_i)\psi  = A^{-1}\left(2(x_i-x_{N+1})\partial_{x_i}-4\chi_i+2\chi_i\right)\psi \\&= A^{-1}\left(2(x_i-x_{N+1})\partial_{x_i}-2\chi_i\right)\Upsilon,  
    \\&(-\xi_i^2\partial_{\xi_i} + \chi_i\xi_i)\psi = A^{-1}(t_i-t_{N+1})\left(\partial_{x_i}-\frac{2\chi_i}{x_i-x_{N+1}} +\frac{2\chi_i}{x_i-x_{N+1}}\right)\eta \\&= A^{-1}(t_i-
    t_{N+1})\partial_{x_i}\eta.
    \end{split}
\end{equation}
Applying this for \(\psi = \Psi\) and \(\eta = \Upsilon\) and substituting into the first and third Ward identities in \eqref{eq:Ward_psi_in_xi} and multiplying with \(A\),
combined with the first Ward identity \eqref{eq:first_Ward_for_Y}, we recover the other two Ward identities \(0 = \sum_{i = 1}^{N+1}t_i\partial_{x_i}\Upsilon\) and
\begin{equation}
    \begin{aligned}
    0 = \sum_{i=1}^{N}(-2x_i \partial_{x_i}+2\chi_i)\Upsilon -2x_{N+1}\partial_{x_{N+1}}\Upsilon =\sum_{i=1}^{N+1}(-2x_i \partial_{x_i}+2\chi_i)\Upsilon - k\Upsilon
    \end{aligned}
\end{equation}
where \(\chi_{N+1}= k/2\).

\subsubsection{Explicit calculation of the KZ equations.}\quad 
We now proceed to reconstruct the KZ equations in the variables $(x_1,\dots,x_{N+1})$ from the KZ equations in the variables $(\xi_1,\dots,\xi_N)$. To do so we apply \eqref{eq:product_rule} and \eqref{eq:lie_action_coordinate_change} twice to see that:
\begin{flalign*}
    &\Tilde{\Omega}_{ij}\Psi = \left(\frac{1}{2}\rho_{\chi_i}^{(i)}(h)\rho_{\chi_j}^{(j)}(h) + \rho_{\chi_i}^{(i)}(e)\rho_{\chi_j}^{(j)}(f) + \rho_{\chi_i}^{(i)}(f)\rho_{\chi_j}^{(j)}(e)\right)\Psi \\
    &= \bigg(\frac{1}{2}(-2\xi_i\partial_{\xi_i}+2\chi_i)(-2\xi_j\partial_{\xi_j}+2\chi_j) 
    + \partial_{\xi_i}(-\xi_j^2\partial_{\xi_j}+2\chi_j\xi_j)+(-\xi_i^2\partial_{\xi_i}+2\chi_i\xi_i)\partial_{\xi_j}\bigg)\Psi \\\nonumber
    &= A^{-1}\bigg(4((x_i-x_{N+1})\partial_{x_i}-\chi_i)((x_j-x_{N+1})\partial_{x_j}-\chi_j)\\\nonumber 
    &+ \frac{(x_i-x_{N+1})^2}{(t_i-t_{N+1})}\left(\partial_{x_i}-\frac{2\chi_i}{(x_i-x_{N+1})}\right)(t_j-t_{N+1}) \partial_{x_j} \\\nonumber
    &+((t_i-t_{N+1})\partial_{x_i})\frac{(x_j-x_{N+1})^2}{(t_j-t_{N+1})}\left(\partial_{x_j}-\frac{2\chi_j}{x_j-x_{N+1}}\right)\bigg)\Upsilon.
\end{flalign*}
This can be written as
\begin{equation*}
    \begin{split}
        &A^{-1}\left(B + \frac{t_j-t_{N+1}}{t_i-t_{N+1}}C + \frac{t_i-t_{N+1}}{t_j-t_{N+1}}D\right)\Upsilon
        \\&= A^{-1} \left(B + C + D + \frac{t_j-t_i}{t_i-t_{N+1}}C + \frac{t_i-t_j}{t_j-t_{N+1}}D\right)\Upsilon,
    \end{split}
\end{equation*}
where \(B,C,D\) are the coefficients of the respective rational functions in \((t_1,\dots,t_{N+1})\).
Therefore, we obtain
\begin{equation*}
     \frac{\Tilde{\Omega}_{ij}}{t_i-t_j}\Psi = A^{-1}\left(\frac{B + C + D}{t_i-t_j} - \frac{1}{t_i-t_{N+1}}C + \frac{1}{t_j-t_{N+1}}D\right)\Upsilon.
\end{equation*}
By summing over \(i \in \{1,\dots,N\}\setminus\{j\}\) and subtracting $(k+2)\partial_{t_{N+1}}\Psi$  we can recover the \((N+1)\)-point KZ equations in the variables $(x_1,\dots,x_{N+1})$.

We begin by showing that $B+C+D$ equates to $\Omega_{ij}$:
\begin{flalign*}
    &B + C + D = (2(x_i-x_{N+1})(x_j-x_{N+1})-(x_i-x_{N+1})^2-(x_j-x_{N+1})^2)\partial_{x_i}\partial_{x_j}\\&+\bigg(2\chi_j\frac{(x_i-x_{N+1})^2}{x_j-x_{N+1}} + 2(x_i-x_{N+1})\bigg(-\chi_j-2\chi_j\frac{(x_i-x_{N+1})^2}{x_j-x_{N+1}} 
    +(x_j-x_{N+1})2\chi_j\bigg)\bigg)\partial_{x_i}\\
    &+ \bigg(2\chi_i(x_i-x_{N+1}) + 2\chi_i\frac{(x_j-x_{N+1})^2}{x_i-x_{N+1}}-2\chi_i\frac{(x_j-x_{N+1})^2}{x_i-x_{N+1}}+ 2(x_j-x_{N+1})(-3\chi_i)\bigg)\partial_{x_j} \\
    & - 6\chi_i\chi_j-4\chi_i\chi_j\frac{x_j-x_{N+1}}{x_i-x_{N+1}}+\chi_i\chi_j\frac{x_j-x_{N+1}}{x_i-x_{N+1}}.
\end{flalign*}
Therefore, we indeed recover
$$\Omega_{ij}= 2(x_i-x_j)\partial_{x_i}\partial_{x_j} + 2(x_i-x_j)(\chi_j\partial_{x_j}-\chi_i\partial_{x_i})+2\chi_i\chi_j.$$

Recalling that 
\begin{align}
    &C = (-(x_i-x_{N+1})^2\partial_{x_i}+ 2\chi_i(x_i-x_{N+1}))\partial_{x_j}\\
    &D =  \partial_{x_i} (-(x_j-x_{N+1})^2\partial_{x_j}+ 2\chi_j(x_j-x_{N+1}))
\end{align}
holds, we continue by calculating
\begin{align*}
    \sum^{N+1}_{\substack{j=1\\ i \neq j}} \frac{C\Upsilon}{t_i-t_{N+1}} &= \frac{(-(x_i-x_{N+1})^2\partial_{x_i}+2\chi_i(x_i-x_{N+1}))}{t_i-t_{N+1}}\sum^n_{\substack{j=0\\ i \neq j}} \partial_{x_j}\Upsilon &&\\
    & = \frac{(x_i-x_{N+1})^2\partial_{x_i}-2\chi_i(x_i-x_{N+1}))}{t_i-t_{N+1}}(\partial_{x_i}+ \partial_{x_{N+1}})\Upsilon
\end{align*}
here we used that 
$$\sum^{N+1}_{\substack{j=1\\ i \neq j}}\partial_{x_j}\Upsilon = -(\partial_{x_i}+ \partial_{x_{N+1}})
\Upsilon.$$ 
Similarly, we calculate
\begin{align*}
    \sum^n_{\substack{j=0\\ i \neq j}} \frac{D}{t_j-t_{N+1}}\Upsilon &=  \partial_{x_i} \sum^n_{\substack{j=0\\ i \neq j}} \frac{-(x_j-x_{N+1})^2\partial_{x_j}+ 2\chi_j(x_j-x_{N+1})}{t_j-t_{N+1}}\Upsilon &&\\
    & = \frac{\partial_{x_i}(-(x_i-x_{N+1})^2\partial_{x_i}+2\chi_i(x_i-x_{N+1}))}{t_i-t_{N+1}}\Upsilon,
\end{align*}
where we employed the first Ward identity. By combing these terms we consider:
\begin{align*}
    \mathcal{E}\Upsilon &\coloneqq \sum^{N+1}_{\substack{j=1\\ i \neq j}} \left(\frac{-C}{t_i-t_{N+1}} + \frac{D}{t_j-t_{N+1}} \right)\Upsilon  
    \\&= \frac{-(x_i-x_{N+1})^2\partial_{x_i}+2\chi_i(x_i-x_{N+1}))}{t_i-t_{N+1}}(\partial_{x_i}+ \partial_{x_{N+1}})\Upsilon \\
    & +\frac{\partial_{x_i}(-(x_i-x_{N+1})^2\partial_{x_i}+2\chi_i(x_i-x_{N+1}))}{t_i-t_{N+1}}\Upsilon \\ 
    & = \frac{1}{t_{N+1}-t_i}\bigg((x_i-x_{N+1})^2\partial_{x_i}^2 - 2s_i(x_i-x_{N+1})\partial_{x_i}  \\
    &+ (x_i-x_{N+1})^2\partial_{x_i}\partial_{x_{N+1}} - 2\chi_i(x_i-x_{N+1})\partial_{x_{N+1}}\\ 
    & - 2(x_i-x_{N+1})\partial_{x_i}+2\chi_i -(x_i-x_{N+1})^2\partial_{x_i}^2 + 2\chi_i(x_i-x_{N+1})\partial_{x_i}\bigg)\Upsilon 
    \\ & = -\frac{1}{t_i-t_{N+1}}\bigg(-(x_i-x_{N+1})^2\partial_{x_i}\partial_{x_{N+1}} \\
    &+2\chi_i(x_i-x_{N+1})\partial_{x_{N+1}}+2(x_i-x_{N+1})\partial_{x_i}-2\chi_i)\bigg)\Upsilon.
\end{align*}

Now let us first recall that:
\begin{equation}
    (k+2)\partial_{t_i}\Psi = \sum^N_{\substack{j=1\\ i \neq j}} \frac{\Tilde{\Omega}_{ij}}{t_i-t_j}\Psi = A^{-1}\sum^n_{\substack{j=0\\ i \neq j}} \frac{\Omega_{ij}}{t_i-t_j}+ A^{-1}\mathcal{E}\Upsilon.
    \label{eq:3}
\end{equation}
On the other hand we have
\begin{equation}
    (k+2)\partial_{t_i} \Upsilon = A(k+2)\partial_{t_i}\Psi- \chi_i\frac{k+2}{t_i-t_{N+1}} \Upsilon - A\frac{k+2}{x_i-x_{N+1}}\partial_{\xi_i}\Psi 
    \label{eq:4}
\end{equation}
and if we now substitute \eqref{eq:3} into \eqref{eq:4} we get the following expression: 
\begin{flalign*}
    (k+2)\partial_{t_i} \Upsilon & = \sum^N_{\substack{j=1\\ i \neq j}} \frac{\Omega_{ij}}{t_i-t_j}\Upsilon + \left(\mathcal{E} + \chi_i\frac{k+2}{t_i-t_{N+1}}\right) \Upsilon -  A\frac{k+2}{x_i-x_{N+1}}\partial_{\xi_i}\Psi \\
    & = \sum^N_{\substack{j=1\\ i \neq j}} \frac{\Omega_{ij}}{t_i-t_j} + \Bigg(\mathcal{E} + \chi_i\frac{k+2}{t_i-t_{N+1}}  - (k+2) \frac{x_i-x_{N+1}}{t_i-t_{N+1}}\left(\partial_{x_i}-\frac{2\chi_i}{x_i-x_{N+1}}\right)\Bigg)\Upsilon \\
    & = \sum^{N+1}_{\substack{j=1\\ i \neq j}} \frac{\Omega_{ij}}{t_i-t_j}\Upsilon + \frac{1}{t_i-t_{N+1}}\bigg(-(x_i-x_{N+1})^2 \partial_{x_i}\partial_{x_{N+1}} \\
    & + 2\chi_i(x_i-x_{N+1})\partial_{x_{N+1}}-2k(x_i-x_{N+1})\partial_{x_i} +k\chi_i\bigg)\Upsilon.
\end{flalign*}
Therefore, we have now fully recovered the \((N+1)\)-point KZ equations for $\Upsilon$ where the highest weight $k/2$ at the \((N+1)\)th point provides the correct $\Omega_{i0}$:
\begin{equation}
    (k+2)\partial_{t_i} \Upsilon = \sum^n_{\substack{j=0\\ i \neq j}} \frac{\Omega_{ij}}{t_i-t_j}\Upsilon + \frac{\Omega_{i0}}{t_i-t_{N+1}}\Upsilon.
\end{equation}

\appendix
\section{Notation}

\subsection{Algebro-geometric}
\begin{itemize}
    \item \(\cO = \bC[\![t]\!] \subset \cK = \bC(\!(t)\!)\) are the rings of formal power series and formal Laurent series respectively and \(D = \textnormal{Spec}(\cO) \supset \textnormal{Spec}(\cK) = D^\circ\) is the formal disc and formal disc without origin respectively. 

    \item \(X\) is a smooth irreducible projective curve. For any finite subset \(\underline{z} = (z_1,\dots,z_N) \in X^N = X \times \dots \times X\) (\(N\)-fold product). Then \(D_{\underline{z}} = \bigsqcup_{i = 1}^N D_{z_i}\) is the formal neighborhood of \(\{z_1,\dots,z_N\}\) in \(X\), \(X^\circ_{\underline{z}} = X \setminus \{z_1,\dots,z_n\}\), and \(D^\circ_{\underline{z}} = D_{\underline{z}} \setminus \{z_1,\dots,z_N\}\).

    \item \(\cO_S\) is the sheaf of regular functions on a scheme \(S\) and for \(\underline{z} \in X^N\) we denote by \(\cO_{\underline{z}} = \bigoplus_{i = 1}^N \cO_{z_i}\) and \(\cK_{\underline{z}} = \bigoplus_{i = 1}^N\cK_{z_i}\) the sections \(\cO_X(D_{\underline{z}})\) and \(\cO_X(D_{\underline{z}}^\circ)\) respectively.

    \item For any ind-schemes \(T\) and \(S\), let us denote \(T(S)\) as the space of morphisms \(S \to T\). If \(S=  \textnormal{Spec}(R)\) is an affine scheme, we also write \(T(S) = T(R)\). In particular, if additionally \(T\) is the scheme associated to a finite-dimensional complex vector space, we have
    \begin{equation}
        T(S) = T(R) = T \otimes R.    
    \end{equation}

    \item The fibered product of two schemes \(S_1,S_2\) with respect to a diagram \(S_1 \to T \leftarrow S_2\) is denoted by \(S_1 \times_T S_2\) and if e.g.\ \(S_1\) is viewed as a scheme over \(T\) in the context, we also use the base change notation \(S_1|_{S_2} = S_1 \times_T \times S_2\) and more specifically if \(S_2 = \{p\}\) is a point, we write \(S_1|_{p} \coloneqq S_1|_{\{p\}}\). If \(H\) is an algebraic group acting on \(S_1\) and \(S_2\) from the left, we also denote the balanced product 
    \(S_1 \times_H S_2 \coloneqq H\setminus(S_1 \times S_2)\). In this paper, these two products can be distinguished by the fact that we never take a fibered product over an algebraic group.
    If we have a diagram \(S_1 \to T \leftarrow S_2\), where an algebraic group \(H\) acts on \(S_1\) and \(S_2\) from the left and acts trivially on \(T\), we write \(S_1 \times_{T,H} S_2 \coloneqq H\setminus(S_1\times_T S_2)\) for the balanced and fibered product.

    \item \(X^\sim\) is the \(\textnormal{Aut}(\cO)\)-bundle over \(X\) consisting of pairs \((x,s)\) of a point \(x \in X\) with the formal coordinate \(s\). For a \(G\)-bundle \(P \to X\), \(P^\sim \to X^\sim\) is the \(\lpgrpp\)-bundle consisting of triples \(((x,s),\sigma)\) of \((x,s) \in X^\sim\) and \(\sigma \in \Gamma(D_x,P)\).
\end{itemize}

\subsection{Lie theoretic}
 \begin{itemize}
     \item \(G\) is a connected simple algebraic group and \(\fg\) is the Lie algebra of \(G\).

    \item We fix Cartan and Borel subalgebras \(\fh \subset \fb_+ = \fh \oplus \fn_+ \subset \fg\) and denote by \(\fg = \fh \oplus \bigoplus_{\alpha \in \Phi} \fg_\alpha\) the associated root space decomposition for the associated polarized root system \(\Phi = \Phi_+ \cup \Phi_-\).

    \item \({}^L\Lambda \subseteq\{\lambda \in \fh \mid \alpha(\lambda) \in \mathbb{Z} \textnormal{ for all }\alpha \in \Phi_+\}\) is the set of all coweights of \(G\) and \({}^L\Lambda_{\ge 0}\) is the subset of dominant coweights, i.e.\ those \(\lambda \in {}^L\Lambda\) satisfying \(\alpha(\lambda)\ge 0\) for all \(\alpha \in \Phi_+\). For any \(\lambda \in {}^L\Lambda\), we denote the associated element of \(G(\textnormal{Spec}(\bC[t,t^{-1}]))\) by \(t^\lambda\).
 
     \item \(\mathfrak{Vir}\) and \(\widehat{\fg}\) denote the Virasoro algebra and the affine Kac-Moody algebra, i.e.\ the canonical central extensions of \(\textnormal{Der}(\cK)\) and \(\fg(\cK) = \fg(D^\circ) =  \fg \otimes \cK\) respectively. We write \(L_m = -t^{m+1}\partial_t \in \mathfrak{Vir}\),  we chose bases \(\{I_j\}_{j = 1}^d \subseteq \fg\), \(\{I^{j}\}_{j = 1}^d \subseteq \fg^*\)  that are dual to each other, and wrote \(a_m \coloneqq a \otimes t^m \in \widehat{\fg}\) for any \(a \in \fg\).

    \item For any subgroup \(K \subseteq G(\cO)\), \(V_{K,k}\) denotes the associated vacuum module of \(\widehat{\fg}\) of level \(k \in \bC\) defined in \eqref{eq:general_vacuum}.

     \item \(U_k(\widehat{\fg})\) is the completed universal enveloping algebra of \(\widehat{\fg}\) at level \(k \in \bC\) and \(V_{\lpgrpp,k} \subset U_k(\widehat{\fg})\) is the vacuum vertex algebra of level \(k\). From Section \ref{sec:coinvariants} onwards, \(k \neq -h^\vee\) and we denote the image of \(L_m\in \mathfrak{Vir}\) in \(V_{\lpgrpp,k}\) given by the Segal-Sugawara construction by \(S_m\). 

     \item For a \(\widehat{\fg}\)-module \(M\) and \(g \in \lpgrp\), \(M^g\) is the \(g\)-twisted module, which has the same underlying base space but is equipped with the twisted action \(x\cdot m \coloneqq (\textnormal{Ad}(g^{-1})x)m\). If \(g = t^\lambda\) for a coweight \(\lambda\), we write \(M^g = M^\lambda\). In particular, the \(\lambda\)-twisted vacuum module is written as \(V_{G(\cO),k}^\lambda\).
\end{itemize}    

\subsection{Moduli spaces}
\begin{itemize}
    \item \(T^\times X\) is the bundle of non-vanishing tangent vectors of \(X\) and 
    \begin{equation}
        C_N^X \coloneqq \{\{(z_i,v_i)\in  T^\times X^N\mid z_i \neq z_j, i,j \in \{1,\dots,N\}\}    
    \end{equation}
    is the configuration space of \(N\)-points of \(X\). For every family \(\underline{K} = (K_1,\dots,K_N)\) of nice (see Section \ref{sec:additional_structures} for a definition) subgroups \(K_1,\dots,K_N \subseteq \lpgrpp\) and any \(G\)-bundle \(P\), let \(C^P_{\underline{K}}\) be the preimage of \(C^X_N\) under the canonical projection \(\prod_{i = 1}^N P^\sim/K_i \to T^\times X^N\).
    
    \item For a \(G\)-bundle \(P\), a family \(\underline{K} = (K_1,\dots,K_N)\) of nice subgroups of \(\lpgrpp\), and a family \(\underline{M} = (M_1,\dots,M_N)\) of \(\underline{K}\)-admissible \(\widehat{\fg}\)-modules (see Section \ref{sec:vbundles_hatg_modules_II}),
    \begin{equation}
        \Delta^{P,\sim}_{\underline{K}}(\underline{M}) \coloneqq (\Delta^{P,\sim}_{{K}_1}(M_1) \boxtimes \dots \boxtimes \Delta^{P,\sim}_{K_N}(M_N))|_{C_{\underline{K}}^P}
    \end{equation} 
    is the associated bundle over \(C_{\underline{K}}^P\), where 
    \begin{equation}
        \Delta_{K_i}^{P,\sim}(M_i) \coloneqq P^\sim \times_{\textnormal{Aut}_+(\cO) \ltimes K_i} M_i.
    \end{equation}
    Furthermore, \[\Delta^P_{\underline{K}}(\underline{M}) = \Delta^{P,\sim}_{\underline{K}}(\underline{M})/\widehat{\fg}^P_{\textnormal{out},N}\cdot\Delta^{P,\sim}_{\underline{K}}(\underline{M}) \] 
    is the sheaf of coinvariants and \(\mathcal{C}^P_{\underline{K}}(\underline{M})\) its dual sheaf of conformal blocks.

    \item \(\textnormal{Bun}_G^X\) is the moduli space of \(G\)-bundles on \(X\). For a family \(\underline{K} = (K_1,\dots,K_N)\) of nice subgroups of \(\lpgrpp\), \(\textnormal{Bun}_{G,\underline{K}}^X\) is the moduli space of \(G\)-bundles with \(\underline{K}\)-structure at the \(N\) marked points. Furthermore, \(\textnormal{Bun}_{G,N}^{X,\sim}\) is the moduli space of \(G\)-bundles with formal trivializations at \(N\) marked points which come equipped with formal coordinates.

    \item For a family \(\underline{K} = (K_1,\dots,K_N)\) of nice subgroups of \(\lpgrpp\) and a family \(\underline{M} = (M_1,\dots,M_N)\) of \(\underline{K}\)-admissible \(\widehat{\fg}\)-modules,  
    \begin{equation}
        \Delta^\sim_{\underline{K}}(\underline{M}) = \textnormal{Bun}_{G,N}^{X,\sim} \times_{\prod_{i = 1}^N(\textnormal{Aut}_+(\cO) \ltimes K_i)} (M_1 \otimes \dots \otimes M_N)
    \end{equation} 
    is the associated bundle over \(\textnormal{Bun}_{G,\underline{K}}^{X}\). Furthermore, 
    \[\Delta_{\underline{K}}(\underline{M}) = \Delta^{\sim}_{\underline{K}}(\underline{M})/\widehat{\fg}_{\textnormal{out},N}\cdot\Delta^{\sim}_{\underline{K}}(\underline{M}) \] 
    is the sheaf of coinvariants and \(\mathcal{C}_{\underline{K}}(\underline{M})\) its dual sheaf of conformal blocks.
\end{itemize}

\newpage
 \printbibliography

\end{document}